%% file: 2002-5.tex
\documentclass{gtart}

\input gtspec

\lognumber{214}
\volumenumber{6}\papernumber{5}\volumeyear{2002}
\pagenumbers{91}{152}
\received{20 October 2001}
\accepted{20 February 2002}
\published{14 March 2002}
\proposed{Walter Neumann}
\seconded{Shigeyuki Morita, Robion Kirby}

\usepackage[all]{xy}
\usepackage{amsmath}
\usepackage{amsthm}
\usepackage{amscd}
\usepackage{amssymb}

\numberwithin{equation}{section}

\newcommand{\nb}[1]{#1\nobreakdash-}

\newcommand\ds\displaystyle

\theoremstyle{definition}
\newtheorem*{definition}{Definition}

\theoremstyle{plain}
\newtheorem{theorem}{Theorem}[section]
\newtheorem*{theorem*}{Theorem}
\newtheorem{proposition}[theorem]{Proposition}
\newtheorem{lemma}[theorem]{Lemma}
\newtheorem{claim}[theorem]{Claim}
\newtheorem{fact}[theorem]{Fact}
\newtheorem{question}[theorem]{Question}

\DeclareMathOperator{\Out}{Out}
\DeclareMathOperator{\Length}{len}
\DeclareMathOperator{\Homeo}{Homeo}
\DeclareMathOperator{\Isom}{Isom}
\DeclareMathOperator{\Aut}{Aut}
\DeclareMathOperator{\Stab}{Stab}
\DeclareMathOperator\Max{Max}
\DeclareMathOperator\diam{diam}
\DeclareMathOperator\image{image}
\DeclareMathOperator\kernel{kernel}
\DeclareMathOperator\interior{int}
\DeclareMathOperator\Hull{Hull}
\DeclareMathOperator\WHull{WH}
\DeclareMathOperator\Comm{Comm}
\DeclareMathOperator\dsup{d_{\mathrm{sup}}}
\DeclareMathOperator\QD{QD}
\DeclareMathOperator\QI{QI}

\newcommand\R{{\mathbf R}}
\newcommand\reals{\R}
\newcommand\hyp{\mathbf{H}}
\newcommand\complex{{\mathbf C}}
\newcommand\Z{{\mathbf Z}}
\newcommand\solv{{\textsc{solv}}}
\renewcommand\H{{\mathcal{H}}}

\newcommand\inject{\hookrightarrow}
\newcommand\homeo{\approx}
\newcommand\infinity{\infty}
\newcommand{\bdy}{\partial}
\newcommand{\from}{\colon}
\newcommand\composed{\circ}
\newcommand\suchthat{\bigm|}
\newcommand\inv{{-1}}
\newcommand\union{\cup}
\newcommand\absvalue[1]{\left| #1 \right|}
\newcommand\abs[1]{\absvalue{#1}}
\newcommand\norm[1]{\left\| #1 \right\|}
\newcommand\wt\widetilde

\newcommand\Id{\text{Id}}
\newcommand\A{\mathcal A}
\newcommand\B{\mathcal B}
\newcommand\intersect{\cap}
\newcommand\FP{{\mathcal{FP}}}
\renewcommand\P{{\mathbf P}}

\newcommand\restrict{\bigm|}
\newcommand\subgroup{<}
\newcommand\semidirect{\rtimes}
\newcommand\normal{\vartriangleleft}
\newcommand\Teichmuller{Teichm\"uller}
\newcommand\Poincare{Poincar\'e}
\newcommand\MCG{{MCG}}

\newcommand\cross{\times}
\newcommand\ext{{\rm ext}}
\newcommand\Haus{{\mathrm{Haus}}}
\newcommand\F{{\cal F}}
\newcommand\G{{\cal G}}
\renewcommand\O{{\cal O}}
\newcommand\M{{\cal M}}
\newcommand\Mod{{\cal M}}
\newcommand\PMF{\P\MF}
\newcommand\MF{{\cal MF}}
\newcommand\C{\mathcal C}
\renewcommand\S{{\cal S}}
\newcommand\T{{\cal T}}
\newcommand\<\langle
\renewcommand\>\rangle

\newcommand\ray[2]{\overrightarrow{[#1,#2)}}
\newcommand\geodesic[2]{\overleftrightarrow{(#1,#2)}}
\newcommand\tree{\frak t}

\begin{document}

\title{Convex cocompact subgroups of\\mapping class groups}

\author{Benson Farb\\Lee Mosher}
\address{Department of Mathematics, University of Chicago \\
5734 University Ave,
Chicago, Il 60637, USA}
\email{farb@math.uchicago.edu}

\secondaddress{Department of Mathematics and Computer Science \\
Rutgers University,
Newark, NJ 07102, USA}
\secondemail{mosher@andromeda.rutgers.edu}

\asciiaddress{Department of Mathematics, University of Chicago\\
5734 University Ave, Chicago, Il 60637, USA\\and\\
Department of Mathematics and Computer Science\\
Rutgers University, Newark, NJ 07102, USA}

\asciiemail{farb@math.uchicago.edu, mosher@andromeda.rutgers.edu}

\begin{abstract} 
We develop a theory of convex cocompact subgroups of the mapping class
group $\MCG$ of a closed, oriented surface $S$ of genus at least $2$,
in terms of the action on \Teichmuller\ space. Given a subgroup $G$ of
$\MCG$ defining an extension $1 \to \pi_1(S) \to\Gamma_G \to G \to 1$,
we prove that if $\Gamma_G$ is a word hyperbolic group then $G$ is a
convex cocompact subgroup of $\MCG$. When $G$ is free and convex
cocompact, called a \emph{Schottky subgroup} of $\MCG$, the converse
is true as well; a semidirect product of $\pi_1(S)$ by a free group
$G$ is therefore word hyperbolic if and only if $G$ is a Schottky
subgroup of $\MCG$. The special case when $G=\Z$ follows from
Thurston's hyperbolization theorem. Schottky subgroups exist in
abundance: sufficiently high powers of any independent set of
pseudo-Anosov mapping classes freely generate a Schottky subgroup.
\end{abstract}

\asciiabstract{ 
We develop a theory of convex cocompact subgroups of the mapping class
group MCG of a closed, oriented surface S of genus at least 2, in
terms of the action on Teichmuller space. Given a subgroup G of MCG
defining an extension L_G:

1---> pi_1(S) ---> L_G ---> G --->1

we prove that if L_G is a word hyperbolic group then G is a convex
cocompact subgroup of MCG. When G is free and convex cocompact, called
a "Schottky subgroup" of MCG, the converse is true as well; a
semidirect product of pi_1(S) by a free group G is therefore word
hyperbolic if and only if G is a Schottky subgroup of MCG. The special
case when G=Z follows from Thurston's hyperbolization theorem.

Schottky subgroups exist in abundance: sufficiently high powers of any
independent set of pseudo-Anosov mapping classes freely generate a
Schottky subgroup.
}


\maketitlepage

\section{Introduction}
\label{SectionIntro}

\subsection{Convex cocompact groups}
\label{SectionConvexCompact}

A \emph{convex cocompact} subgroup of $\Isom(\hyp^n)$, the isometry group of hyperbolic
$n$--space, is a discrete subgroup $G \subgroup \Isom(\hyp^n)$, with limit set $\Lambda_G
\subset\bdy\hyp^n$, such that $G$ acts cocompactly on the convex hull $\Hull_G \subset
\hyp^n$ of its limit set $\Lambda_G$. It follows that $G$ is a word hyperbolic group with
model geometry $\Hull_G$ and Gromov boundary $\Lambda_G$. Given any
finitely generated, discrete subgroup $G\subgroup\Isom(\hyp^n)$, $G$ is
convex cocompact if and only if any orbit of $G$ is a quasiconvex subset
of $\hyp^n$. Convex cocompact subgroups satisfy several useful
properties: every infinite order element of $G$ is loxodromic;
$\Lambda_G$ is the smallest nontrivial $G$--invariant closed subset of
$\overline\hyp^n = \hyp^n \union \bdy
\hyp^n$; the action of $G$ on $\bdy\hyp^n\setminus\Lambda_G$ is properly discontinuous; assuming
$\Lambda_G \ne \bdy\hyp^n$, the stabilizer subgroup of $\Lambda_G$ is a finite index supergroup
of $G$, and it is the relative commensurator of $G$ in $\Isom(\hyp^n)$.

A \emph{Schottky group} is a convex cocompact subgroup of $\Isom(\hyp^n)$ which is free.
Schottky subgroups of $\Isom(\hyp^n$) exist in abundance and can be constructed using the
classical ping-pong argument, attributed to Klein: if $\phi_1,\ldots,\phi_n$ are loxodromic
elements whose axes have pairwise disjoint endpoints at infinity, then sufficiently high
powers of $\phi_1,\ldots,\phi_n$ freely generate a Schottky group.\footnote{The term
``Schottky group'' sometimes refers explicitly to a subgroup of $\Isom(\hyp^n)$ produced by
the ping-pong argument, but the broader reference to free, convex cocompact subgroups has
become common.}

We shall investigate the notions of convex cocompact groups and Schottky groups in the context of
\Teichmuller\ space. Given a closed, oriented surface $S$ of genus $\ge 2$, the mapping class
group $\MCG$ acts as the full isometry group of the \Teichmuller\ space $\T$
\cite{Royden:IsomTeichReport}.\footnote{In this paper, $\MCG$ includes orientation
reversing mapping classes, and so represents what is sometimes called the ``extended''
mapping class group.} This action extends to the Thurston compactification
$\overline\T = \T \union \PMF$ \cite{FLP}. \Teichmuller\ space is \emph{not} Gromov hyperbolic
\cite{MasurWolf:TeichNotHyp}, no matter what finite covolume, equivariant metric one picks
\cite{BrockFarb:curvature}, and yet it exhibits many aspects of a hyperbolic metric space
\cite{Minsky:quasiprojections} \cite{MasurMinsky:complex1}. A general theory of limit
sets of finitely generated subgroups of $\MCG$ is developed in \cite{McCarthyPapa:dynamics}.

In this paper we develop a theory of convex cocompact subgroups and Schottky subgroups of
$\MCG$ acting on $\T$, and we show that Schottky subgroups exist in abundance. We apply this
theory to relate convex cocompactness of subgroups of $\MCG$ with the large scale
geometry of extensions of surface groups by subgroups of $\MCG$.

Our first result establishes the concept of convex cocompactness for subgroups of $\MCG$, by
proving the equivalence of several properties:

\begin{theorem}[Characterizing convex cocompactness]
\label{TheoremQuasiconvex}
Given a finitely generated subgroup $G\subgroup\MCG$, the following statements are
equivalent:
\begin{itemize}
\item Some orbit of $G$ is quasiconvex in $\T$.
\item Every orbit of $G$ is quasiconvex in $\T$.
\item $G$ is word hyperbolic, and there is a $G$--equivariant embedding 
$\bdy f \from \bdy G \to \PMF$ with image $\Lambda_G$ such that the following properties hold:
\begin{itemize} 
\item Any two distinct points $\xi,\eta \in \Lambda_G$ are the ideal endpoints of a unique
geodesic $\geodesic{\xi}{\eta}$ in $\T$.
\item Let $\WHull_G$ be the ``weak hull'' of $G$, namely the union of the geodesics
$\geodesic{\xi}{\eta}$, $\xi \ne \eta \in \Lambda_G$. Then the action of $G$ on $\WHull_G$ is
cocompact, and if $f \from G \to \WHull_G$ is any $G$--equivariant map then $f$ is a
quasi-isometry and the following map is continuous:
$$\bar f =  f \union \bdy f \from G \union \bdy G \to \overline\T = \T \union \PMF
$$ 
\end{itemize}
\end{itemize} 
\end{theorem}
Any such subgroup $G$ is said to be \emph{convex cocompact}. This theorem is proved in
Section~\ref{SectionEquivalence}.

A convex cocompact subgroup $G \subgroup \MCG$ shares many properties with convex cocompact
subgroups of $\Isom(\hyp^n)$. Every infinite order element of $G$ is pseudo-Anosov
(Proposition~\ref{PropCnvxCbddPsA}). The limit set $\Lambda_G$ is the smallest nontrivial closed
subset of $\overline T$ invariant under the action of $G$, and the action of
$G$ on $\PMF-\Lambda_G$ is properly discontinuous (Proposition~\ref{PropCnvxCbddLimitSet}); this
depends on work of McCarthy and Papadoupolos \cite{McCarthyPapa:dynamics}. The stabilizer of
$\Lambda_G$ is a finite index supergroup of $G$ in $\MCG$, and it is the relative commensurator
of $G$ in $\MCG$ (Corollary~\ref{PropConvCbddComm}).

A \emph{Schottky subgroup} of $\MCG = \Isom(\T)$ is defined to be a convex cocompact subgroup
which is free of finite rank. In Theorem~\ref{TheoremSchottkyAbundance} we prove that if
$\phi_1,\ldots,\phi_n$ are pseudo-Anosov elements of $\MCG$ whose axes have pairwise disjoint
endpoints in $\PMF$, then for all sufficiently large positive integers $a_1,\ldots,a_n$ the
mapping classes $\phi_1^{a_1},\ldots,\phi_n^{a_n}$ freely generate a Schottky subgroup of $\MCG$.

\paragraph{Warning} Our formulation of convex cocompactness in $\T$ is not as strong as in
$\hyp^n$. Although there is a general theory of limit sets of finitely generated
subgroups of $\MCG$ \cite{McCarthyPapa:dynamics}, we have no general theory of their convex
hulls. Such a theory would be tricky, and unnecessary for our purposes. In particular, when $G$
is convex cocompact, we do not know whether there is a closed, convex, $G$--equivariant subset of
$\T$ on which $G$ acts cocompactly. One could attempt to construct such a subset by adding to
$\WHull_G$ any geodesics with endpoints in $\WHull_G$, then adding geodesics with endpoints in
that set, etc, continuing transfinitely by adding geodesics and taking closures until the result
stabilizes; however, there is no guarantee that $G$ acts cocompactly on the result.

\subsection{Surface group extensions}
\label{SectionExtensions}

There is a natural isomorphism of short exact sequences
$$
\xymatrix{
1 \ar[r] & \pi_1(S,p) \ar[r]^(.4){\iota} \ar@{=}[d] & \MCG(S,p) \ar[r]^q \ar@2{~>}[d] & \MCG(S)
\ar[r] \ar@2{~>}[d] & 1
\\ 1 \ar[r] & \pi_1(S,p) \ar[r]        & \Aut(\pi_1(S,p)) \ar[r]   & \Out(\pi_1(S,p)) \ar[r] & 1
}$$ where $\MCG(S,p)$ is the mapping class group of $S$ punctured at the
base point $p$. In the bottom sequence, the inclusion $\pi_1(S,p)$ is
obtained by identifying $\pi_1(S,p)$ with its group of inner
automorphisms, an injection since $\pi_1(S,p)$ is centerless. For each
based loop $\ell$ in $S$, $\iota(\ell)$ is the punctured mapping class
which ``pushes'' the base point $p$ around the loop $\ell$ (see
Section~\ref{SectionTeich} for the exact definition). The homomorphism
$q$ is the map which ``forgets'' the puncture $p$. Exactness of the top
sequence is proved in
\cite{Birman:Braids}. The isomorphism $\MCG(S) \approx\Out(\pi_1(S,p))$ follows from work of
Dehn--Nielsen
\cite{Nielsen:untersuchungen}, Baer \cite{Baer:Isotopien}, and Epstein \cite{Epstein:curves}. As a
consequence, either of the above sequences is natural for extensions of $\pi_1(S)$, in the
following sense. For any group homomorphism $G \to\MCG(S)$, by applying the fiber product
construction to the homomorphisms 
$$\xymatrix{
\MCG(S,p) \ar[dr] & & G \ar[dl] \\
 & \MCG(S)
}$$
we obtain a group $\Gamma_G$ and a commutative diagram of short exact sequences
$$
\xymatrix{
1 \ar[r] & \pi_1(S) \ar[r] \ar[d] & \Gamma_G \ar[r] \ar[d] & G \ar[r]
\ar[d] & 1 \\ 1 \ar[r] & \pi_1(S) \ar[r] & \MCG(S,p) \ar[r] & \MCG(S)
\ar[r] & 1 }
$$
Note that we are suppressing the homomorphism $G \to \MCG(S)$ in the notation $\Gamma_G$. If $G$
is free then the top sequence splits and we can write $\Gamma_G = \pi_1(S) \semidirect G$,
where again our notation suppresses a lift $G \to \Aut(\pi_1(S))$ of the given homomorphism $G
\to \MCG(S) \approx \Out(\pi_1(S))$. 

Every group extension $1 \to \pi_1(S) \to E \to G \to 1$ arises from the above
construction, because the given extension determines a homomorphism $G \to \Out(\pi_1(S)) \approx
\MCG(S)$ which in turn determines an extension $1 \to \pi_1(S) \to \Gamma_G
\to G \to 1$ isomorphic to the given extension. 

When $P$ is a cyclic subgroup of $\MCG$, Thurston's hyperbolization
theorem for mapping tori (see, eg,\ \cite{Otal:fibered}) shows 
that $\pi_1(S) \semidirect P$ is the fundamental group of a closed,
hyperbolic 3--manifold if and only if $P$ is a pseudo-Anosov subgroup. In
particular, $\pi_1(S) \semidirect P$ is a word hyperbolic group if and
only if $P$ is a convex cocompact subgroup of $\MCG$. Our results about
the extension groups $\Gamma_G$ are aimed towards generalizing this
statement as much as possible. The theme of these results is that the
geometry of $\Gamma_G$ is encoded in the geometry of the action of $G$
on~$\T$.

From \cite{Mosher:HypExt} it follows that if $\Gamma_G$ is word
hyperbolic then $G$ is word hyperbolic. Our next result gives much more
precise information:

\begin{theorem}[Hyperbolic extension has convex cocompact quotient] 
\label{TheoremHypExtQuotient}
If $\Gamma_G$ is word hyperbolic then the homomorphism $G \to \MCG$ has
finite kernel and convex cocompact image.
\end{theorem} 

This theorem is proved in Section~\ref{SectionExtensionQuotient}. Finiteness of the kernel $K$ of
$G\to \MCG$ is easy to prove, using the fact that $\pi_1(S) \cross K$ is a subgroup of~$\Gamma_G$.
If $K$ is infinite, then either it is a torsion group, or it has an infinite order element and so
$\Gamma_G$ has a $\Z\oplus\Z$ subgroup; in either case, $\Gamma_G$ cannot be word hyperbolic.
Because one can mod out by a finite kernel without affecting word hyperbolicity of the extension
group, this brings into focus the extensions defined by inclusion of subgroups of
$\MCG$.

We are particularly interested in free subgroups of $\MCG$. A finite rank, free, convex
cocompact subgroup is called a \emph{Schottky subgroup}. For Schottky subgroups we have a converse
to Theorem~\ref{TheoremHypExtQuotient}, giving a complete characterization of word hyperbolic
groups $\Gamma_F$ when $F \subgroup \MCG$ is free:

\begin{theorem}[Surface-by-Schottky group has hyperbolic extension]
\label{TheoremSchottkyHypExt}
If $F$ is a finite rank, free subgroup of $\MCG$ then the extension group $\Gamma_F=\pi_1(S)
\semidirect F$ is word hyperbolic if and only if $F$ is a Schottky group.
\end{theorem}

This is proved in Section~\ref{SectionSchottky}. Some special cases of
this theorem are immediate. It is not hard to see that
$\pi_1(S)\semidirect F$ has a $\Z\oplus\Z$ subgroup if and only if there
exists a nontrivial element $f \in F$ which is not pseudo-Anosov. Such
an element $f$, being infinite order, must be reducible. Assuming $f \in
F$ is nontrivial and reducible, the group $\pi_1(S)\semidirect F$
contains the subgroup $\pi_1(S) \semidirect \<f\>$ which is the
fundamental group of a closed \nb{3}manifold that contains an
incompressible torus. Conversely, when $\pi_1(S)\semidirect F$ has a
$\Z\oplus\Z$ subgroup then that subgroup must map onto an infinite
cyclic subgroup $\<f\>\subset F$ whose action on $\pi_1(S)$ preserves a
nontrivial conjugacy class, and so $f$ is not pseudo-Anosov. Theorem
\ref{TheoremSchottkyHypExt} is therefore mainly about free,
pseudo-Anosov subgroups of $\MCG$ (see
Question~\ref{QuestionAllSchottky} below).

The abundance of word hyperbolic extensions of the form $\pi_1(S)
\semidirect F$ was proved in
\cite{Mosher:hypbyhyp}. It was shown by McCarthy \cite{McCarthy:Tits} and Ivanov
\cite{Ivanov:subgroups} that if $\phi_1,\ldots,\phi_n$ are pseudo-Anosov
elements of $\MCG$ which are pairwise independent, meaning that their
axes have distinct endpoints in the Thurston boundary $\PMF$, then
sufficiently high powers of these elements freely generate a
pseudo-Anosov subgroup $F$. The main result of \cite{Mosher:hypbyhyp}
shows in addition that, after possibly making the powers higher, the
group $\pi_1(S) \semidirect F$ is word hyperbolic. The nature of the
free subgroups $F\subgroup \MCG$ produced in
\cite{Mosher:hypbyhyp} was somewhat mysterious, but Theorems~\ref{TheoremHypExtQuotient}
and~\ref{TheoremSchottkyHypExt} clear up this mystery by characterizing the subgroups $F$
using an intrinsic property, namely convex cocompactness. 

By combining \cite{Mosher:hypbyhyp} and Theorem~\ref{TheoremSchottkyHypExt}, we immediately have
the following result:

\begin{theorem}[Abundance of Schottky subgroups]
\label{TheoremSchottkyAbundance}
\quad If $\phi_1,\ldots,\phi_n \in \MCG$ are pairwise independent
pseudo-Anosov elements, then for all sufficiently large positive
integers $a_1,\ldots,a_n$ the mapping classes
$\phi_1^{a_1},\ldots,\phi_n^{a_n}$ freely generate a Schottky subgroup
$F$ of $\MCG$. 
\end{theorem}

Finally, we shall show in Section~\ref{SectionOrbifolds} that all of the
above results generalize to the setting of closed hyperbolic
\nb{2}orbifolds. These generalized results find application in the
results of \cite{FarbMosher:sbf}, as we now recall.

\subsection{An application}
\label{SectionApplication}
In the paper \cite{FarbMosher:sbf} we apply our theory of Schottky subgroups of
$\MCG$ to investigate the large-scale geometry of word hyperbolic surface-by-free groups:

\begin{theorem*}{\rm\cite{FarbMosher:sbf}}\qua
Let $F \subset \MCG(S)$ be Schottky. Then the group $\Gamma_F = \pi_1(S) \semidirect F$ is
quasi-isometrically rigid in the strongest sense: 
\begin{itemize}
\item $\Gamma_F$ embeds with finite index in
its quasi-isometry group $\QI(\Gamma_F)$. 
\end{itemize}
It follows that:
\begin{itemize}
\item Let $H$ be any finitely generated group.  If 
$H$ is quasi-isometric to $\Gamma_F$, then there exists a finite normal
subgroup $N \normal H$ such that $H/N$ and $\Gamma_F$ are abstractly
commensurable.
\item The abstract commensurator group $\Comm(\Gamma_F)$ is isomorphic to $\QI(\Gamma_F)$, and
can be computed explicitly.
\end{itemize}
\end{theorem*}

The computation of $\Comm(\Gamma_F) \approx \QI(\Gamma_F)$ goes as
follows. Among all orbifold subcovers $S \to \O$ there exists a unique
minimal such subcover such that the subgroup $F \subgroup \MCG(S)$
descends isomorphically to a subgroup $F'\subgroup\MCG(\O)$.  The whole
theory of Schottky groups extends to general closed hyperbolic
orbifolds, as we show in Section~\ref{SectionOrbifolds} of this
paper. In particular, $F'$ is a Schottky subgroup of $\MCG(\O)$. By
Corollary~\ref{PropConvCbddComm} it follows that $F'$ has finite index
in its relative commensurator $N \subgroup\MCG(\O)$, which can be
regarded as a virtual Schottky group. The inclusion $N \subgroup
\MCG(\O)$ determines a canonical extension $1 \to \pi_1(\O) \to \Gamma_N
\to N \to 1$, and we show in \cite{FarbMosher:sbf} that the extension
group $\Gamma_N$ is isomorphic to $\QI(\Gamma_F)$.

\subsection{Some questions} 
\label{SectionQuestions}
Our results on convex cocompact and Schottky subgroups of $\MCG$
motivate several questions.

Proposition~\ref{PropCnvxCbddPsA} implies that if $F$ is a Schottky subgroup of $\MCG$ then
every nontrivial element of $F$ is pseudo-Anosov.
\begin{question}
\label{QuestionAllSchottky}
Suppose $F \subgroup \MCG$ is a finite rank, free subgroup all of whose
nontrivial elements are pseudo-Anosov. Is $F$ convex cocompact? In other
words, is $F$ a Schottky group?
\end{question}

A non-Schottky example $F$ would be very interesting for the following
reasons. There exist examples of infinite, finitely presented groups
which are not word hyperbolic and whose solvable subgroups are all
virtually cyclic, but all known examples fail to be of finite type; see for example
\cite{Brady:branchedcoverings}. If there were a non-Schottky subgroup $F\subgroup\MCG$ as in
Question~\ref{QuestionAllSchottky}, then the group $\pi_1(S)\semidirect F$ would be of finite type
(being the fundamental group of a compact aspherical \nb{3}complex), it would not be word
hyperbolic (since $F$ is not Schottky), and every nontrivial solvable subgroup
$H\subgroup\pi_1(S)\semidirect F$ would be infinite cyclic. To see why the latter holds, since
$\pi_1(S)\semidirect F$ is a torsion free subgroup of $\MCG(S,p)$ it follows by
\cite{BirmanLubotzkyMcCarthy} that the subgroup $H$ is finite rank free abelian. Under the
homomorphism
$H\to F$, the groups $\image(H \to F)\subgroup F$ and $\kernel(H \to F)\subgroup \pi_1(S)$ each
are free abelian of rank at most~1, and so it suffices to rule out the case where the image and
kernel both have rank~$1$. But in that case we would have a pseudo-Anosov element of
$\MCG(S)$ which fixes the conjugacy class of some infinite order element of $\pi_1 S$, a
contradiction.


Note that Question~\ref{QuestionAllSchottky} has an analogue in the theory of Kleinian groups: if
$G$ is a discrete, cocompact subgroup of $\Isom(\hyp^3)$, is every free subgroup of $G$ a
Schottky subgroup? More generally, if $G$ is a discrete, cofinite volume subgroup of
$\Isom(\hyp^3)$, is every free loxodromic subgroup of $G$ a Schottky group? The first question,
at least, would follow from Simon's tame ends conjecture \cite{Canary:CoveringSurvey}.

For a source of free, pseudo-Anosov subgroups on which to test
question~\ref{QuestionAllSchottky}, consider Whittlesey's group \cite{Whittlesey:pAgroup}, an
infinite rank, free, normal, pseudo-Anosov subgroup of the mapping class group of a closed,
oriented surface of genus~2.
\begin{question}
Is every finitely generated subgroup of Whittlesey's group a Schottky group?
\end{question}

Concerning non-free subgroups of $\MCG$, note first that Question~\ref{QuestionAllSchottky} can
also be formulated for any finitely generated subgroup of $\MCG$, though we have no examples of
non-free pseudo-Anosov subgroups. This invites comparison with the situation in $\Isom(\hyp^n)$
where it is known for any $n \ge 2$ that there exist convex cocompact subgroups which are not
Schottky, indeed are not virtually Schottky.

\begin{question}
\label{QuestionNonFreeConvexCocompact}
Does there exist a convex cocompact subgroup $G \subgroup \MCG$ which is not
Schottky, nor is virtually Schottky? 
\end{question}

The converse to Theorem~\ref{TheoremHypExtQuotient}, while proved for free
subgroups in Theorem~\ref{TheoremSchottkyHypExt}, remains open in general. This issue becomes
particularly interesting if Question~\ref{QuestionNonFreeConvexCocompact} is answered
affirmatively:

\begin{question} 
\label{QuestionCCImpliesHypExt}
If $G \subgroup \MCG$ is convex cocompact, is the extension group $\Gamma_{G}$ word
hyperbolic? 
\end{question}

Surface subgroups of mapping class groups are interesting. Gonzalez-D\`{\i}ez
and Harvey showed that $\MCG$ can contain the fundamental group of a closed, oriented surface
of genus $\ge 2$ \cite{GonzalezDiezHarvey:surfacesubgroup}, but their construction always
produces subgroups containing mapping classes that are not pseudo-Anosov. 

If questions~\ref{QuestionNonFreeConvexCocompact} and~\ref{QuestionCCImpliesHypExt} were
true, it would raise the stakes on the fascinating question of whether there exist
surface-by-surface word hyperbolic groups:

\begin{question} 
Does there exist a convex cocompact subgroup $G \subgroup \MCG$ isomorphic to the
fundamental group of a closed, oriented surface $S_g$ of genus $g \ge 2$? If so, is the
surface-by-surface extension group $\Gamma_G$ word hyperbolic?
\end{question}

Misha Kapovich shows in \cite{MKapovich:ComplexSurfaces} that when $G$ is a surface
group, the extension group $\Gamma_G$ cannot be a lattice in $\Isom(\complex\hyp^2)$.

\subsection{Sketches of proofs}
\label{SectionSketches}

Although \Teichmuller\ space $\T$ is not hyperbolic in any reasonable sense
\cite{MasurWolf:TeichNotHyp}, \cite{BrockFarb:curvature}, nevertheless it possesses interesting
and useful hyperbolicity properties. To formulate these, recall that the action of $\MCG$ by
isometries on $\T$ is smooth and properly discontinuous, with quotient orbifold $\Mod = \T / \MCG$
called the \emph{moduli space} of $S$. The action is \emph{not} cocompact, and we define a subset
$A\subset\T$ to be \emph{cobounded} if its image under the universal covering map $\T \to \Mod$
has compact closure in $\Mod$, equivalently there is a compact subset of
$\T$ whose translates under $\Isom(\T)$ cover $A$.

In \cite{Minsky:quasiprojections}, Minsky proves (see
Theorem~\ref{TheoremMinskyContraction} below) that if $\ell$ is a
cobounded geodesic in $\T$ then any projection $\T \to \ell$ that takes
each point of $\T$ to a closest point on $\ell$ satisfies properties
similar to a closest point projection from a $\delta$--hyperbolic metric
space onto a bi-infinite geodesic.  This projection property is a key
step in the proof of the Masur--Minsky theorem
\cite{MasurMinsky:complex1} 
that Harvey's curve complex is a $\delta$--hyperbolic metric space.  
These results say intuitively that $\T$
exhibits hyperbolicity as long as one focusses only on cobounded
aspects. Keeping this in mind, the tools of
\cite{Minsky:quasiprojections} and
\cite{MasurMinsky:complex1} can be used to prove Theorem~\ref{TheoremQuasiconvex} along the
classical lines of the proof for subgroups of $\Isom(\hyp^n)$.

The proof of Theorem~\ref{TheoremSchottkyHypExt}, that $\pi_1(S) \semidirect F$ is word
hyperbolic if $F$ is Schottky, uses the Bestvina--Feighn combination theorem
\cite{BestvinaFeighn:combination}. Consider a tree $\tree$ on which $F$
acts freely and cocompactly, and choose an $F$--equivariant mapping $\phi\from
\tree\to\T$. Let $\H \to \T$ be the canonical hyperbolic plane bundle over \Teichmuller\
space. Pulling back via $\phi$ we obtain a hyperbolic plane bundle $\pi\from\H_\tree \to
\tree$, and $\pi_1(S) \semidirect F$ acts properly discontinuously and cocompactly on
$\H_\tree$. This shows that $\H_\tree$ is a model geometry for the group $\pi_1(S)
\semidirect F$, and in particular $\H_\tree$ is a $\delta$--hyperbolic metric space if and only
if $\pi_1(S) \semidirect F$ is word hyperbolic.

By the Bestvina--Feighn combination theorem \cite{BestvinaFeighn:combination} and its converse
due to Gersten \cite{Gersten:cohomological}, hyperbolicity of $\H_\tree$ is equivalent to
$\delta$--hyperbolicity of each ``hyperplane'' $\H_\ell = \pi^\inv(\ell)$, where $\ell$ ranges
over all the bi-infinite lines in $\tree$ and $\delta$ is independent of $\ell$.

Recall that for each \Teichmuller\ geodesic $g$, the canonical marked Riemann surface
bundle $\S_g$ over $g$ carries a natural \emph{singular \solv\ metric}; the bundle $\S_g$
equipped with this metric is denoted $\S^\solv_g$. Lifting the metric to the universal cover
$\H_g$ we obtain a singular \solv\ space denoted $\H^\solv_g$.

When $F$ is a Schottky group, convex cocompactness tells us that for each bi-infinite geodesic
$\ell$ in $\tree$, the map $\ell\xrightarrow{\phi}\T$ is a quasigeodesic and there is a
unique \Teichmuller\ geodesic $g$ within finite Hausdorff distance from $\phi(\ell)$. This
feeds into Proposition~\ref{PropFellowTravellers}, a basic construction principle for
quasi-isometries which will be used several times in the paper. The conclusion is:

\begin{fact}
The hyperplane $\H_\ell$ is uniformly quasi-isometric to the singular \solv--space
$\H^\solv_g$, by a quasi-isometry which is a lift of a closest point map $\ell \to g$. 
\end{fact}

Uniform hyperbolicity of singular \solv--spaces $\H^\solv_g$, where $g$ is a uniformly
cobounded geodesic in $\T$, is then easily checked by another application of the
Bestvina--Feighn combination theorem, and Theorem~\ref{TheoremSchottkyHypExt} follows.

For Theorem~\ref{TheoremHypExtQuotient}, we first outline the proof in the special case of a
free subgroup of $\MCG$. As noted above, using Gersten's converse to the Bestvina--Feighn
combination theorem, word hyperbolicity of $\pi_1(S) \semidirect F$ implies uniform
hyperbolicity of the hyperplanes $\H_\ell$. Now we use a result of Mosher
\cite{Mosher:StableQuasigeodesics}, which shows that from uniform hyperbolicity of the
hyperplanes $\H_\ell$ it follows that the lines $\ell$ are uniform quasigeodesics in $\T$,
and each $\ell$ has uniformly finite Hausdorff distance from some \Teichmuller\ geodesic $g$.
Piecing together the geodesics $g$ in $\T$, one for each geodesic $\ell$ in $\tree$, we obtain
the data we need to prove that $F$ is Schottky.

The general proof of Theorem~\ref{TheoremHypExtQuotient} follows the same outline, except
that we cannot apply Gersten's converse to the Bestvina--Feighn combination theorem. That
result applies only to the setting of groups acting on trees, not to the setting of
Theorem~\ref{TheoremHypExtQuotient} where $\Gamma_G$ acts on the Cayley graph of $G$. To
handle this problem we need a new idea: a generalization of Gersten's converse to the
Bestvina--Feighn combination theorem, which holds in a much broader setting. This
generalization is contained in Lemma~\ref{LemmaHypImpliesFlare}. The basis of this
result is an analogy between the ``flaring property'' of Bestvina--Feighn and the divergence
of geodesics in a word hyperbolic group \cite{Cannon:TheoryHyp}.

\paragraph{Acknowledgements} We are grateful to the referee for a thorough
reading of the paper, and for making numerous useful comments. 

Both authors are supported in part by the National Science Foundation.

\section{Background}
\label{SectionPreliminaries}

\subsection{Coarse language}
\label{SectionCoarseLanguage}

\paragraph{Quasi-isometries and uniformly proper maps}

Given a metric space $X$ and two subsets $A,B \subset X$, the 
\emph{Hausdorff distance}
$d_\Haus(A,B)$ is the infimum of all real numbers $r$ such that each point of $A$ is within
distance $r$ of a point of $B$, and vice versa.

A \emph{quasi-isometric embedding} between two metric spaces $X,Y$ is a
map $f \from X \to Y$ such that for some $K \ge 1$, $C \ge 0$, we have
$$\frac{1}{K} d(x,y) - C \le d(fx,fy) \le K d(x,y) + C $$ for each $x,y
\in X$. To refer to the constants we say that $f$ is a
$K,C$--quasi-isometric embedding.

For example, a quasigeodesic embedding $\R\to X$ is
called a
\emph{quasigeodesic line} in $X$. We also speak of \emph{quasigeodesic rays or
segments} with the domain is a half-line or a finite segment,
respectively. Since every map of a segment is a quasi-isometry, it
usually behooves one to include the constants and speak about a
$(K,C)$--quasi-isometric segment.

A \emph{quasi-isometry} between two metric spaces $X,Y$ is a map $f
\from X \to Y$ which, for some $K \ge 1$, $C \ge 0$ is a $K,C$
quasi-isometry and has the property that $\image(f)$ has Hausdorff
distance $\le C$ from $Y$. Every quasi-isometry $f\from X \to Y$ has a
\emph{coarse inverse}, which is a quasi-isometry $\bar f \from Y \to X$
such that $\bar f \composed f \from X
\to X$ is a bounded distance in the sup norm from $\Id_X$, and similarly for $f \composed \bar f
\from Y \to Y$; the sup norm bounds and the quasi-isometry constants of $\bar f$ depend only on
the quasi-isometry constants of $f$.

More general than a quasi-isometric embedding is a \emph{uniformly proper embedding} $f \from X
\to Y$, which means that there exists $K \ge 1$, $C \ge 0$, and a function $r \from [0,\infinity)
\to [0,\infinity)$ satisfying $r(t) \to \infinity$ as $t \to \infinity$, such that
$$r(d(x,y)) \le d(fx,fy) \le K d(x,y) + C
$$
for each $x,y \in X$.

\paragraph{Geodesic and quasigeodesic metric spaces}
A metric space is \emph{proper} if closed balls are compact. A metric $d$ on a space $X$ is
called a \emph{path metric} if for any $x,y \in X$ the distance $d(x,y)$ is the infimum of
the path lengths of rectifiable paths between $x$ and $y$, and $d$ is called a \emph{geodesic
metric} if $d(x,y)$ equals the length of some rectifiable path between $x$ and $y$. The
following fact is an immediate consequence of the Ascoli--Arzela theorem:

\begin{fact}
\label{FactGeodesicMetricSpace} A compact path metric space is a geodesic metric space. More
generally, a proper path metric is a geodesic metric.
\qed\end{fact}

The Ascoli--Arzela theorem also shows that for any proper geodesic metric space $X$, every path
homotopy class contains a shortest path. This implies that the metric on $X$ lifts to a
geodesic metric on any covering space of $X$. 

A metric space $X$ is called a \emph{quasigeodesic metric space} if there exists constants
$\lambda,\epsilon$ such that for any $x,y \in X$ there exists an interval $[a,b] \subset \reals$
and a $\lambda,\epsilon$ quasigeodesic embedding $\sigma \from [a,b] \to X$ such that
$\sigma(a)=x$ and $\sigma(b)=y$. For example, if $Y$ is a geodesic metric space and $X$ is a
subset of $Y$ such that $d_\Haus(X,Y)<\infinity$ then $X$ is a quasigeodesic metric space.

The fundamental theorem of geometric group theory, first known to Efremovich, to Schwarzc, and to
Milnor, can be given a general formulation as follows. Let $X$ be a proper, quasigeodesic metric
space, and let the group $G$ act on $X$ properly discontinuously and cocompactly, by an action
denoted $(g,x)\mapsto g \cdot x$. Then $G$ is finitely generated, and for any
base point $x_0 \in X$ the map $G \to X$ defined by $g \mapsto g \cdot x_0$ is a quasi-isometry
between the word metric on $G$ and the metric space $X$.

\paragraph{Uniform families of quasi-isometries} The next lemma says a family of geodesic
metrics which is ``compact'' in a suitable sense has the property that any two metrics in the
family are uniformly quasi-isometric, with respect to the identity map.

Given a compact space $X$, let $M(X)$ denote the space of metrics generating the topology of $X$,
regarded as a subspace of $[0,\infinity)^{X\cross X}$ with the topology of uniform convergence.

\begin{lemma}
\label{LemmaFamilyOfMetrics}
Let $X$ be a compact, path connected space with universal cover $\wt X$. Let $D \subset
M(X)$ be a compact family of geodesic metrics. Let $\wt D$ be the set of lifted metrics on $\wt
X$. Then there exist $K \ge 1$, $C \ge 0$ such that for any $\wt d,\wt d' \in \wt D$ the identity
map on $\wt X$ is a $K,C$ quasi-isometry between $(\wt X,\wt d)$ and $(\wt X,\wt d')$.
\end{lemma}

\begin{proof} By compactness of $D$, the metric spaces $X_d$ have a uniform injectivity
radius---that is, there exists $\epsilon>0$ such that for each $d \in D$ every homotopically
nontrivial closed curve in $X_d$ has length $> 4\epsilon$, and it follows that every closed
$\epsilon$ ball in $X_d$ lifts isometrically to $\wt X_d$. Let $P \subset \wt X \cross \wt X$ be
the set of pairs $(x,y) \in \wt X \cross \wt X$ such that for some $d \in \wt D$ we have $d(x,y)
\le\epsilon$. Evidently $\pi_1(X)$ acts cocompactly on $P$, and so we have a finite
supremum
\begin{align*}
A &= \sup\{\wt d(x,y) \suchthat \wt d \in \wt D \text{ and } (x,y) \in P\}
\end{align*}
Given $\wt d \in \wt D$ and $x,y \in \wt X$, choose a $\wt d$--geodesic $\gamma$ from $x$ to
$y$ and let $x=x_0,x_1,\ldots,x_{n-1},x_n=y$ be a monotonic sequence along $\gamma$
such that $d(x_{i-1},x_i) = \epsilon$ for $i=1,\ldots,n-1$ and $d(x_{n-1},x_n) \le \epsilon$.
For any $\wt d' \in \wt D$ it follows that:
$$\wt d'(x,y) \le An = A \left\lceil \frac{\wt d(x,y)}{\epsilon} \right\rceil \le
\frac{A}{\epsilon} \wt d(x,y) + A
$$
Setting $K=\frac{A}{\epsilon}$ and $C=A$ the lemma follows.
\end{proof}

\paragraph{Hyperbolic metric spaces} A geodesic metric space $X$ is \emph{hyperbolic} if there
exists $\delta \ge 0$ such that for any $x,y,z \in X$ and any geodesics $\overline{xy}$,
$\overline{yz}$, $\overline{zx}$, any point on $\overline{xy}$ has distance $\le\delta$ from some
point on $\overline{yz} \union \overline{zx}$. A finitely generated group is \emph{word
hyperbolic} if the Cayley graph of some (any) finite generating set, equipped with the geodesic
metric making each edge of length~1, is a hyperbolic metric space.

If $X$ is $\delta$--hyperbolic, then for any $\lambda \ge 1$, $\epsilon \ge 0$ there exists $A$,
depending only on $\delta,\lambda,\epsilon$, such that the following hold: for any $x,y \in X$,
any $\lambda,\epsilon$ quasigeodesic segment between $x$ and $y$ has Hausdorff distance $\le A$
from any geodesic segment between $x$ and $y$; for any $x \in X$, any $\lambda,\epsilon$
quasigeodesic ray starting at $x$ has Hausdorff distance $\le A$ from some geodesic ray starting
at $x$; and any $\lambda,\epsilon$ quasigeodesic line in $X$ has Hausdorff distance $\le A$ from
some geodesic line in $X$.

The \emph{boundary} of $X$, denoted $\bdy X$, is the set of coarse equivalence classes of geodesic
rays in $X$, where two rays are coarsely equivalent if they have finite Hausdorff distance. For
any $\xi\in\bdy X$ and $x_0 \in X$, there is a ray based at $x_0$ representing $\xi$; we denote
such a ray $\ray{x_0}{\xi}$. For any $\xi \ne \eta \in \bdy X$ there is a geodesic line $\ell$ in
$X$ such that any point on $\ell$ divides it into two rays, one representing $\xi$ and the other
representing $\eta$. 

Assuming $X$ is proper, there is a compact topology on $X\union \bdy X$ in which $X$ is dense,
which is characterized by the following property: a sequence $\xi_i \in X\union\bdy X$ converges
to $\xi\in \bdy X$ if and only if, for any base point $p \in X$, 
if $\ray{p}{\xi_i}$ denotes
either a segment from $p$ to $\xi_i \in X$, or a ray from $p$ with ideal endpoint $\xi_i \in \bdy
X$, then any subsequential limit of the sequence $\ray{p}{\xi_i}$ is a ray with ideal endpoint
$\xi$. It follows that any quasi-isometric embedding between $\delta$--hyperbolic geodesic metric
spaces extends to a continuous embedding of boundaries. In particular, if $X$ is hyperbolic then
the action of $\Isom(X)$ on $X$ extends continuously to an action on $X \union \bdy X$.

The following fundamental fact is easily proved by considering what happens to geodesics in a
$\delta$--hyperbolic metric space under a quasi-isometry.

\begin{lemma}
\label{LemmaBdyValueUniqueness}
For all $\delta \ge 0$, $K \ge 1$, $C \ge 0$ there exists $A \ge 0$ such
that the following holds. If $X,Y$ are two $\delta$--hyperbolic metric
spaces and if $f, g \from X \to Y$ are two $K,C$ quasi-isometries such
that $\bdy f = \bdy g \from \bdy X \to \bdy Y$, then: 
$$\dsup(f,g) =
\sup_{x \in X} d(f(x),g(x)) \le A\eqno{\qed} $$
\end{lemma}

\subsection{\Teichmuller\ space and the Thurston boundary}
\label{SectionTeich}

Fix once and for all a closed, oriented surface $S$ of genus $g \ge 2$.
Let $\C$ be the set of isotopy classes of nontrivial simple closed
curves on~$S$.

The fundamental notation for the paper is as follows. Let $\T$ be the
\Teichmuller\ space of $S$. Let $\MF$ be the space of measured
foliations on $S$, and let $\PMF$ be the space of projective measured
foliations on $S$, with projectivization map $\P \from \MF \to
\PMF$. The Thurston compactification of \Teichmuller\ space is
$\overline \T = T \union \PMF$. Let $\MCG$ be the mapping class group of
$S$, and let $\Mod = \T / \MCG$ be the moduli space of $S$. Definitions
of these objects are all recalled below.

The \Teichmuller\ space $\T$ is the set of hyperbolic structures on $S$
modulo isotopy, with the structure of a smooth manifold diffeomorphic to
$\R^{6g-6}$ given by Fenchel--Nielsen coordinates.  The Riemann mapping
theorem associates to each conformal structure on $S$ a unique
hyperbolic structure in that conformal class, and hence we may naturally
identify $\T$ with the set of conformal structures on $S$ modulo
isotopy. Given a conformal structure or a hyperbolic structure $\sigma$,
we will often confuse $\sigma$ with its isotopy class by writing
$\sigma\in \T$.

There is a length pairing $\T \cross \C \to \R_+$ which associates to each $\sigma \in \T$, $C \in
\C$ the length of the unique simple closed geodesic on the hyperbolic surface $\sigma$ in the
isotopy class $C$. We obtain a map $\T \to [0,\infinity)^\C$ which is an embedding
with image homeomorphic to an open ball of dimension $6g-6$. Moreover, under projectivization
$[0,\infinity)^\C \to \P[0,\infinity)^\C$, $\T$ embeds in $\P[0,\infinity)^\C$ with precompact
image. 

\paragraph{Thurston's boundary}

A \emph{measured foliation} $\F$ on $S$ is a foliation with finitely many singularities
equipped with a positive transverse Borel measure, with the property that for each singularity $s$
there exists $n \ge 3$ such that in a neighborhood of $s$ the foliation $\F$ is modelled on the
horizontal measured foliation of the quadratic differential $z^{n-2} dz^2$ in the complex plane. A
\emph{saddle connection} of $\F$ is a leaf segment connecting two distinct singularities;
collapsing a saddle connection to a point yields another measured foliation on $S$. The set of
measured foliations on $S$ modulo the equivalence relation generated by isotopy and saddle
collapse is denoted $\MF$. Given a measured foliation $\F$, its equivalence class is denoted $[\F]
\in
\MF$; elements of $\MF$ will often be represented by the letters $X,Y,Z$.

For each measured foliation $\F$, there is a function $\ell_\F \from \C \to
[0,\infinity)$ defined as follows. Given a simple closed curve $c$, we may pull back the
transverse measure on $\F$ to obtain a measure on $c$, and then integrate over $c$ to obtain a
number $\int_c\F$.  
Define $\ell_\F(c)=i(\F,c)$ to be the infimum of $\int_{c'}\F$ as $c'$ ranges over the isotopy
class of $c$. The function $\ell_\F$ is well-defined up to equivalence, thereby defining an
embedding $\MF \to [0,\infinity)^\C$ whose image is homeomorphic to $\R^{6g-6} - \{0\}$. 

Given a measured foliation $\F$, multiplying the transverse measure by a positive
scalar $r$ defines a measured foliation denoted $r \cdot \F$, yielding a positive scalar
multiplication operation $\reals \cross \MF\to \MF$. With respect to the equivalence relation
$\F\sim r \cdot \F$, $r > 0$, the set of equivalence classes is denoted
$\PMF$ and the projection is denoted $\P \from \MF \to \PMF$. We obtain an embedding $\PMF
\to \P[0,\infinity)^\C$ whose image is homeomorphic to a sphere of dimension
$6g-7$. We often use the letters $\xi,\eta,\zeta$ to represent points of $\PMF$.

Thurston's compactification theorem \cite{FLP} says, by embedding into
$\P[0,\infinity)^\C$, that there is a homeomorphism of triples: 
$$(\overline\T, \T, \PMF) \approx (B^{6g-6},\interior(B^{6g-6}),S^{6g-7})
$$ 
We will also need the
standard embedding $\C\to\MF$, defined on $[c]$ as follows. Take an embedded annulus $A \subset S$
foliated by circles in the isotopy class $[c]$, and assign total transverse measure~1 to the
annulus. Choose a deformation retraction of each component of the closure of $S-A$ onto a finite
1--complex, and extend to a map $f \from S \to S$ homotopic to the
identity and which is an embedding
on
$\interior(A)$. The measured foliation on $A$ pushes forward under $f$ to the desired measured
foliation on $S$, giving a well-defined point in $\MF$ depending only on $[c]$.

The intersection number $\MF\cross\C\xrightarrow{i(\cdot,\cdot)} [0,\infinity)$ extends
continuously to $\MF\cross \MF \xrightarrow{i(\cdot,\cdot)} [0,\infinity)$. This intersection
number is most efficaciously defined in terms of measured geodesic laminations. 

\paragraph{Marked surfaces} Having fixed once and for all the surface $S$, a \emph{marked
surface} is a pair $(F,\phi)$ where $F$ is a surface and $\phi \from S \to F$ is a
homeomorphism. Thus we may speak about a marked hyperbolic surface, a marked Riemann surface,
a marked measured foliation on a surface, etc. 

Given a marked hyperbolic surface $\phi \from S \to F$, pulling back via $\phi$ determines a
hyperbolic structure on $S$ and a point of $\tree$. Two marked hyperbolic surfaces $\phi \from S
\to F$ and $\phi' \from S \to F'$ give the same element of $\T$ if and only if they are
equivalent in the following sense: there exists an isometry $h \from F \to F'$ such that
$\phi'{}^\inv\composed h \composed \phi \from S \to S$ is isotopic to the identity. In this
manner, we can identify the collection of marked hyperbolic surfaces up to equivalence with
the \Teichmuller\ space $\T$ of $S$. This allows us the freedom of representing a point of
$\T$ by a hyperbolic structure on some other surface $F$, assuming implicitly that we have a
marking $\phi\from S \to F$. The same discussion holds for marked Riemann surfaces, marked
measured foliations on surfaces, etc.

Given two marked surfaces $\phi \from S \to F$, $\phi' \from S \to F'$, a
\emph{marked map} is a homeomorphism $\psi\from F \to F'$ such that $\psi \composed
\phi$ is isotopic to $\phi'$.

\paragraph{Mapping class groups and moduli space}

Let $\Homeo(S)$ be the group of homeomorphisms of $S$, let $\Homeo_0(S)$ be the normal
subgroup consisting of homeomorphisms isotopic to the identity, and let $\MCG =
\MCG(S) = \Homeo(S) / \Homeo_0(S)$ be the \emph{mapping class group} of $S$. Pushing a
hyperbolic structure on $S$ forward via an element of $\Homeo(S)$ gives a well-defined action of
$\MCG$ on $\T$. This action is smooth and properly discontinuous but not cocompact. It
follows that the \emph{moduli space} $\Mod =\T /\MCG$ is a smooth, noncompact orbifold with
fundamental group $\MCG$ and universal covering space $\T$.

Let $\Homeo(S,p)$ be the group of homeomorphisms of $S$ preserving a base point $p$, let
$\Homeo_0(S,p)$ be the normal subgroup consisting of those homeomorphisms which are isotopic
to the identity leaving $p$ stationary, and let $\MCG(S,p) = \Homeo(S,p) / \Homeo_0(S,p)$.
Recall the short exact sequence: 
$$
1 \to \pi_1(S,p) \xrightarrow{\iota} \MCG(S,p) \xrightarrow{q} \MCG(S) \to 1
$$
The map
$q$ is the map which ``forgets'' the puncture $p$. To define the map $\iota$, for each closed
loop $\ell \from [0,1] \to S$ based at $p$, choose numbers $0=x_0 < x_1 <
\ldots < x_n=1$ and embedded open balls $B_1,\ldots,B_n \subset S$ so that $\ell[x_{i-1},x_i]
\subset B_i$ for $i=1,\ldots,n$, and let $\pi_i \from S \to S$ be a homeomorphism which is the
identity on $S-B_i$ and such that $\pi_i(\ell(x_{i-1})) = \ell(x_i)$. Then $\iota(\ell)$ is
defined to be the isotopy class rel $p$ of the homeomorphism $\pi_n \composed \pi_{n-1} \composed
\cdots\composed \pi_1 \from (S,p) \to (S,p)$, which we say is obtained by ``pushing'' the point
$p$ around the loop~$\ell$. The mapping class $\iota(\ell)$ is well-defined independent of the
choices made, and independent of the choice of $\ell$ in its path homotopy class. When $\ell$ is
simple, $\iota(\ell)$ may also be described as the composition of opposite Dehn twists on the two
boundary components of a regular neighborhood of~$\ell$. For details see \cite{Birman:Braids}. 

As noted in the introduction, by the
Dehn--Nielsen--Baer--Epstein theorem, the above sequence is naturally isomorphic to the sequence 
$$
1 \to \pi_1(S,p) \to \Aut(\pi_1(S,p)) \to \Out(\pi_1(S,p)) \to 1
$$

\paragraph{Canonical bundles} Over the 
\Teichmuller\ space $\T$ of $S$ there is a
\emph{canonical marked hyperbolic surface bundle} $\S \to \T$, defined
as follows.  Topologically $\S = S \cross\T$, with the obvious marking
$S\xrightarrow{\approx} S\cross
\sigma = \S_\sigma$ for each $\sigma\in \T$. As $\sigma$ varies over
$\T$, one can assign a hyperbolic structure on $S$ in the class of
$\sigma$, varying continuously in the $C^\infinity$ topology on
Riemannian metrics; this follows from the description of Fenchel--Nielsen
coordinates. It follows that on each fiber $\S_\sigma$ of $\S$ there is
a hyperbolic structure which varies continuously in $\sigma$. Note that
by the Riemann mapping theorem we can also think of $\S$ as the
canonical marked Riemann surface bundle over $\T$.

The action of $\MCG$ on $\T$ lifts uniquely to an action on $\S$, such
that for each fiber $\S_\sigma$ and each $[h] \in \MCG$ the map
$$\S_\sigma \xrightarrow{[h]} \S_{[h](\sigma)} $$ is an isometry, and
the map $$S \xrightarrow{\approx} \S_\sigma \xrightarrow{[h]}
\S_{[h](\sigma)} \xrightarrow{\approx} S $$ is in the mapping class
$[h]$. 


The universal cover of $\S$ is called the canonical $\hyp^2$--bundle over
$\T$, denoted $\H \to \T$. There is a fibration preserving, isometric
action of the once-punctured mapping class group $\MCG(S,p)$ on the total space $\H$
such that the quotient action of $\MCG(S,p)$ on $\S$ has kernel
$\pi_1(S,p)$, and corresponds to the given action of $\MCG(S) =
\MCG(S,p) / \pi_1(S,p)$ on $\S$. Also, the action of $\pi_1(S,p)$ on any
fiber of $\H$ is conjugate to the action on the universal cover $\wt S$
by deck transformations. Bers proved in \cite{Bers:FiberSpaces}
that $\H$ is a \Teichmuller\ space in its own right: there is an $\MCG(S,p)$ equivariant
homeomorphism between $\H$ and the \Teichmuller\ space of the once-punctured surface $S-p$.

The tangent bundle $T\S$ has a smooth \nb{2}dimensional \emph{vertical sub-bundle} $T_v\S$
consisting of the tangent planes to fibers of the fibration $\S\to\T$. A \emph{connection} on the
bundle $\S\to\T$ is a smooth codimension--2 sub-bundle of $T\S$ complementary to $T_v\S$. The
existence of an $\MCG$--equivariant connection on $\S$ can be derived following standard methods,
as follows. Choose a locally finite, equivariant open cover of $\T$, and an equivariant partition
of unity dominated by this cover. For each $\MCG$--orbit of this cover, choose a representative $U
\subset\T$ of this orbit, and choose a linear retraction $T\S_U \to T_v\S_U$. Pushing these
retractions around by the action of $\MCG$ and taking a linear combination using the partition of
unity, we obtain an equivariant linear retraction $T\S \to T_v\S$, whose kernel is the desired
connection.

By lifting to $\H$ we obtain a connection on the bundle $\H\to\T$, equivariant with respect to
the action of the group $\MCG(S,p)$.

\subparagraph{Notation} Given any subset $A \subset \T$, or more generally any continuous map $A
\to
\T$, by pulling back the bundle $\S\to\T$ we obtain a bundle $\S_A \to A$, as shown in the
following diagram:
$$\xymatrix{
\S_A \ar@{.>}[d] \ar@{.>}[r]  & \S \ar[d] \\
A \ar[r]            & \T
}$$
Similarly, the pullback of the bundle $\H\to\T$ is denoted $\H_A \to A$.

\paragraph{Quadratic differentials} 

Given a conformal structure $\sigma$
on $S$, a \emph{quadratic differential} $q$ on $\S_\sigma$ assigns to each conformal
coordinate~$z$ an expression of the form $q(z) dz^2$ where $q(z)$ is a complex valued function on
the domain of the coordinate system, and
$$q(z) \left(\frac{dz}{dw}\right)^2 = q(w), \quad\text{for overlapping coordinates $z,w$.}
$$
We shall always assume that the functions $q(z)$ are holomorphic, in other words, our quadratic
differentials will always be ``holomorphic'' quadratic differentials. A quadratic differential
$q$ is \emph{trivial} if $q(z)$ is always the zero function.

Given a nontrivial quadratic differential $q$ on $\S_\sigma$, a point $p
\in \S_\sigma$ is a zero of $q$ in one coordinate if and only if it is a
zero in any coordinate; also, the order of the zero is well-defined. If
$p$ is not a zero then there is a coordinate $z$ near~$p$, unique up to
multiplication by $\pm 1$, such that $p$ corresponds to the origin and
such that $q(z) \equiv 1$; this is called a \emph{regular canonical
coordinate}. If $p$ is a zero of order $n \ge 1$ then up to
multiplication by the $(n+2)^{\text{nd}}$ roots of unity there exists a
unique coordinate~$z$ in which $p$ corresponds to the origin and such
that $q(z) = z^n$; this is called a
\emph{singular canonical coordinate}. There is a well-defined \emph{singular Euclidean metric}
$\abs{q(z)}\abs{dz}^2$ on $S$, which in any regular canonical coordinate $z=x+iy$ takes the
form $dx^2 + dy^2$. In any singular canonical coordinate this metric has finite area, and so
the total area of $S$ in this singular Euclidean metric is finite, denoted $\norm{q}$. We say
that $q$ is \emph{normalized} if $\norm{q}=1$. 

By the Riemann--Roch theorem, the quadratic differentials on $\S_\sigma$ form a complex vector
space $\QD_\sigma$ of complex dimension $3g-3$, and these vector spaces fit together, one for
each $\sigma\in\T$, to form a complex vector bundle over $\T$ denoted $\QD \to \T$. \Teichmuller\
space has a complex structure whose cotangent bundle is canonically isomorphic to the bundle
$\QD$. The \Teichmuller\ metric on $\T$ induces a Finsler metric on the (real) tangent bundle of
$\T$, and the norm $\norm{q}$ is dual to this metric. The normalized quadratic differentials form
a sphere bundle $\QD^1\to\T$ of real dimension $6g-7$ embedded in $\QD$. 

Corresponding to each quadratic differential $q$ on $\S_\sigma$ there is a pair of measured
foliations, the \emph{horizontal foliation} $\F_x(q)$ and the \emph{vertical foliation}
$\F_y(q)$. In a regular canonical coordinate $z=x+iy$, the leaves of $\F_x(q)$ are parallel to
the $x$--axis and have transverse measure $\abs{dy}$, and the leaves of $\F_y(q)$ are parallel
to the $y$--axis and have transverse measure $\abs{dx}$. The foliations $\F_x(q)$, $\F_y(q)$
have the zero set of $q$ as their common singularity set, and at each zero of order $n$ both
have an \emph{$(n+2)$--pronged singularity}, locally modelled on the singularity at
the origin of the horizontal and vertical measured foliations of $z^n dz^2$.

Conversely, consider a \emph{transverse pair of measured foliations} $(\F_x,\F_y)$ on $S$
which means that $\F_x, \F_y$ have the same singular set, are transverse at all regular
points, and at each singularity $s$ there is a number $n \ge 3$ such that $\F_x$ and
$\F_y$ are locally modelled on the horizontal and vertical measured foliations of $z^{n-2}
dz^2$. Associated to the pair $\F_x,\F_y$ there are a conformal structure and a quadratic
differential defined as follows. Near each regular point, there is an oriented coordinate
$z=x+iy$ in which $\F_x$ is the horizontal foliation with transverse measure $\abs{dy}$, and
$\F_y$ is the vertical foliation with transverse measure $\abs{dx}$. These regular
coordinates have conformal overlap. Near any singularity $s$, at which $\F_x$, $\F_y$ are
locally modelled on the the horizontal and vertical foliations of $z^n dz^2$, the coordinate
$z$ has conformal overlap with any regular coordinate. We therefore obtain a conformal
structure $\sigma(\F_x,\F_y)$ on $S$, on which we have a quadratic differential $q(\F_x,\F_y)$
defined in regular coordinates by $dz^2$. 

A pair of measured foliations $(X,Y) \in \MF(F) \cross \MF(F)$ is said to \emph{jointly fill} the
surface $F$ if, for every $Z \in \MF(F)$, either $i(X,Z) \ne 0$ or $i(Y,Z) \ne 0$. This condition
is invariant under positive scalar multiplication on $\MF(F)$, and so joint filling is
well-defined for a pair of points in $\PMF(F)$. A basic fact is that a pair $X,Y \in \MF(F)$
jointly fills $F$ if and only if there exist a transverse pair of measured foliations
$\F_x,\F_y$ representing $X,Y$; moreover, such a pair $\F_x,\F_y$ is unique up to \emph{joint
isotopy}, meaning that for any other transverse pair $\F'_x,\F'_y$ representing $X,Y$
respectively, there exists $h\in \Homeo_0(S)$ such that $\F'_x = h(\F_x)$, $\F'_y=h(\F_y)$. These
facts may be proved by passing back and forth between measured geodesic laminations and
measured foliations.

By uniqueness up to joint isotopy as just described, it follows that for each
jointly filling pair $(X,Y) \in \MF(F) \cross \MF(F)$ there is a conformal structure
$\sigma(\F_x,\F_y)$ and quadratic differential $q(\F_x,\F_y)$ on $\sigma(X,Y)$, well-defined
up to isotopy independent of the choice of a transverse pair $\F_x,\F_y$ representing $X,Y$.
We thus have a well-defined point $\sigma(X,Y) \in \T(F)$ and a well-defined element $q(X,Y)
\in \QD_{\sigma(X,Y)} \T(F)$.

\paragraph{Geodesics and a metric on $\T$}

We shall describe geodesic lines in $\T$ following \cite{GardinerMasur}
and \cite{HubbardMasur:qd}; of course everything depends on \Teichmuller's theorem (see eg,\
\cite{Abikoff:realanalytic} or \cite{ImayoshiTaniguchi}).

Let $\FP\subset\MF\cross \MF$ denote the set of jointly filling
pairs, and let $\P\FP$ be the image of $\FP$ under the product of projection maps
$\P\cross\P\from\MF\cross\MF \to \PMF\cross\PMF$. 

Associated to each jointly filling pair $(\xi,\eta) \in \P\FP$ we associate a
\emph{\Teichmuller\ line} $\geodesic{\xi}{\eta}$, following \cite{GardinerMasur}. Choosing a
transverse pair of measured foliations $\F_x,\F_y$ representing $\xi,\eta$ respectively, we
obtain a \emph{parameterized \Teichmuller\ geodesic} given by the map $t \mapsto \sigma(e^{-t}
\F_x, e^{t} \F_y)$; it follows from \Teichmuller's theorem that this map is an embedding
$\reals \to\T$. Uniqueness of $\F_x,\F_y$ up to joint isotopy and positive scalar multiplication
imply that the map $t \mapsto \sigma(e^{-t} \F_x, e^{t} \F_y)$ is well-defined up to translation
of the $t$--parameter, as is easily checked. Thus, the image of this map is well defined and is
denoted $\geodesic{\xi}{\eta}$; in addition, parameter difference between points on the line is
well-defined, and there is a well-defined orientation. The \emph{positive direction} of the
geodesic is defined to be the point $\eta=\P\F_y \in \PMF$, the projective class of the vertical
measured foliation; the negative direction is the point $\xi=\P\F_x\in\PMF$. Note that as $t\to
+\infinity$ the vertical measured foliation becomes ``exponentially thicker'' and so dominates
over the horizontal foliation which becomes ``exponentially thinner'', a useful mnemonic for
remembering which direction is which. 

\Teichmuller's theorem says that any two distinct points of $\T$ lie on a unique
\Teichmuller\ line: for any $\sigma \ne \tau \in \T$ there exists a unique pair $(\xi,\eta)
\in\P\FP$ such that $\sigma,\tau \in \geodesic{\xi}{\eta}$. Moreover, if $d(\sigma,\tau)$ is the
parameter difference between $\sigma$ and $\tau$ along this geodesic, then $d$ is a metric on
$\T$, called the \emph{\Teichmuller\ metric}. In particular, each line $\geodesic{\xi}{\eta}$ is,
indeed, a geodesic for the \Teichmuller\ metric. It is also true that the segment
$[\sigma,\tau]\subset\geodesic{\xi}{\eta}$ is the unique geodesic segment connecting $\sigma$ to
$\tau$, and hence geodesic segments are uniquely extensible. Thus we obtain a 1--1 correspondence
between oriented geodesic segments and the set $\T \cross \T$. Also, every bi-infinite geodesic
line in $\T$ is uniquely expressible in the form $\geodesic{\xi}{\eta}$, and so we obtain a 1--1
correspondence between oriented geodesic lines and the set $\P\FP \subset \PMF \cross \PMF$.

There is a also 1--1 correspondence between geodesic rays in $\T$ and the set $\T \cross
\PMF$: for any $\sigma \in\T$ and $\eta\in\PMF$ there is a unique geodesic ray, denoted
$\ray{\sigma}{\eta}$, whose endpoint is $\sigma$ and whose direction is $\eta \in
\PMF$, and every geodesic ray has this form. This is an immediate consequence of the
Hubbard--Masur theorem \cite{HubbardMasur:qd}, which says that for each $\sigma \in \T$ the map
$\QD_\sigma \to \MF$ taking $q \ne 0\in\QD_\sigma$ to $[\F_y(q)]$ is a homeomorphism.

Throughout the paper, the term ``geodesic'' will refer to any geodesic
segment, ray, or line in $\T$. Geodesics in $\T$ are \emph{uniquely
extendable}: any geodesic segment or ray is contained in a unique
geodesic line. Since $\T$ is a complete metric space, an argument using
the Ascoli--Arzela theorem shows that any sequence of geodesics, each
element of which intersects a given bounded subset of $\T$, has a
subsequence converging pointwise to a geodesic.

By unique extendability of geodesics it follows that $\T$ is a proper, geodesic metric space.
From the definitions it follows that the action of $\MCG$ on $\T$ is isometric, and so the
metric on $\T$ descends to a proper, geodesic metric on $\Mod = \T / \MCG$.

The reader is cautioned that a geodesic ray $\ray{\sigma}{\eta}$ is \emph{not known} to
converge in $\overline\T$ to its direction $\eta\in\PMF$. However, consider the case where $\eta$
is \emph{uniquely ergodic}, which means that for any measured foliation $\F$ representing $\eta$,
every transverse measure on the underlying singular foliation of $\F$ is a scalar multiple of the
given measure on $\F$. In this case the ray $\ray{\sigma}{\eta}$ does converge to $\eta$, as
is proved by Masur \cite{Masur:TwoBoundaries}, and so in this situation the direction $\eta$ is
also called the \emph{end} or \emph{endpoint} of the ray.

\paragraph{Cobounded geodesics in $\T$} 

A subset $A \subset \T$ is \emph{cobounded} if the image of $A$ under the projection $\T \to
\Mod$ is a bounded subset of $\Mod$; equivalently, there is a bounded subset of $\T$ whose
translates by the action of $\MCG$ cover $A$. If the bounded set $\B\subset\Mod$ contains
the projected image of $A$ then we also say that $A$ is \emph{$\B$--cobounded}. Since $\Mod$ is
a proper metric space it follows that $A$ is cobounded in $\T$ if and only if $A$ is
``co-precompact'', meaning that the projection of $A$ to $\Mod$ has compact closure.

One common gauge for coboundedness, as noted by Mumford \cite{Mumford:Mahler}, is the injectivity
radius of a hyperbolic structure, or to put it another way, the length $\ell(\sigma)$ of the
shortest closed geodesic in a hyperbolic structure~$\sigma$.\footnote{Also called the ``systole''
in the differential geometry literature.} For each $\epsilon>0$ the ``$\epsilon$--thick subset''
of $\T$, namely the set $\T_\epsilon=\{\sigma\in\T\suchthat\ell(\sigma)\ge\epsilon\}$, is an
$\MCG$ equivariant subset of $\T$ projecting to a compact subset of $\Mod$, and as $\epsilon \to
0$ this gives an exhaustion of $\Mod$ by compact sets. A subset of $\T$ is therefore cobounded if
and only if it is contained in the $\epsilon$--thick subset of $\T$ for some $\epsilon>0$.

Extremal length, rather than hyperbolic length, is used to obtain another common gauge of
coboundedness, and is comparable to the length of the shortest geodesic by Maskit's work
\cite{Maskit:comparison}.

We rarely use any particular gauge for coboundedness. Instead, the primary way in which we use
coboundedness is in carrying out compactness arguments over closed, bounded subsets. For this
reason we rarely refer to any gauge, instead sticking with coboundedness as the more primitive
mathematical concept.

One important fact we need is that if $\rho=\ray{\sigma}{\eta}$ is a cobounded geodesic ray in
\Teichmuller\ space then $\rho$ converges to $\eta$ in Thurston's compactification $\overline \T
= \T \union\PMF$. This follows from two theorems of Masur. 
  First, since $\rho$ is cobounded, the direction $\eta \in \PMF$ is
uniquely ergodic; this result, proved in \cite{Masur:exchange}, was later sharpened
in \cite{Masur:NonergodicDimension} to show that if $\eta$ is not uniquely ergodic
then the projection of $\ray{\sigma}{\eta}$ to moduli space leaves every compact subset. Second,
when $\eta$ is uniquely ergodic, any ray with direction $\eta$ converges to $\eta$ in Thurston's
compactification. This is a small part of a Masur's Two Boundaries Theorem
\cite{Masur:TwoBoundaries}, concerning relations between
the \Teichmuller\ boundary and the Thurston boundary of $\T$ (we will use the full power of this
theorem in the proof of Theorem~\ref{TheoremQuasiconvex}).

The following result is essentially a consequence of \cite{Minsky:quasiprojections}:

\begin{lemma}[End Uniqueness]
\label{LemmaEndUniqueness} 
If $\ray{\sigma}{\xi}$, $\ray{\tau}{\eta}$ are two cobounded rays in
$\T$ which have finite Hausdorff distance in $\T$ then $\xi=\eta$. If $\geodesic{\xi}{\xi'}$,
$\geodesic{\eta}{\eta'}$ are two cobounded lines in $\T$ which have finite Hausdorff distance
then, up to relabelling the ends of one of the lines, we have $\xi=\eta$ and
$\xi'=\eta'$, and so $\geodesic{\xi}{\xi'} = \geodesic{\eta}{\eta'}$.
\end{lemma}

\begin{proof} For the proof we review briefly notions of extremal length, in the classical
setting of simple closed curves, as well as Kerckhoff's extension to the setting of
measured foliations \cite{Kerckhoff:asymptotic}.

Recall that for any conformal structure on an open annulus $A$ there is a unique Euclidean
annulus of the form $S^1 \cross (0,M)$ conformally equivalent to $A$, with $M \in \R_+ \union
\{\infinity\}$; the modulus of $A$, denoted $M(A)$, is defined to be the number~$M$. For any
Riemann surface $\S_\sigma$ and any isotopy class of simple closed curves $[c] \in \C$, the
\emph{extremal length} $\ell_\ext(\sigma,[c])$ is the infimum of $1/M(A)$ taken over all
annuli $A\subset F$ whose core is in the isotopy class $[c]$.

Kerckoff proved \cite{Kerckhoff:asymptotic} that the function $\ell_\ext \from \T
\cross (\R_+ \cdot \C) \to (0,\infinity)$ defined by $\ell_\ext(\sigma,r[c])\mapsto r
\ell_\ext(\sigma,[c])$ extends continuously to a function $\ell_\ext \from \T \cross
\MF\to[0,\infinity)$. Moreover, for any transverse pair of measured foliations $\F_x,\F_y$ with
associated conformal structure $\sigma=\sigma(\F_x,\F_y)$ and quadratic differential
$q=q(\F_x,\F_y)$, we have 
$$\ell_\ext(\sigma,\F_y) = \sqrt{\norm{q}} 
$$


Given $X \in \MF$, the \emph{extremal length horoball} based at $X$ is
defined to be $H(X) = \{\sigma \in \T \suchthat \ell_\ext(\sigma,X) \le
1\}$. Note for example that, setting $\xi = \P X$, for every $\eta \in
\PMF$ the extremal length of $X$ at points of $\geodesic{\eta}{\xi}$
decreases strictly monotonically to zero as the point moves towards
$\xi$, and so every \Teichmuller\ geodesic with positive direction $\P X$
eventually enters $H(X)$ in the positive direction and, once in, never
leaves. Given $\xi \in \PMF$, there is a one parameter family of
extremal length horoballs based at $\xi$, namely $H(X)$ for all $X \in
\MF$ such that $\P X=\xi$.

For the first sentence of the theorem, consider two geodesic rays
$\ray{\sigma}{\xi}$, $\ray{\tau}{\eta}$ such that $\xi \ne \eta \in
\PMF$. Pick any extremal length horoball $H$ based at $\eta$. The proof
of Theorem 4.3 of 
\cite{Minsky:quasiprojections} shows that $H \intersect \ray{\sigma}{\xi}$ is bounded. However, $H
\intersect \ray{\tau}{\eta}$ is an infinite subray of $\ray{\tau}{\eta}$, and moreover as a point
$p \in\ray{\tau}{\eta}$ travels to infinity in $\ray{\tau}{\eta}$ the
horoball $H$ contains a larger and larger ball in $\T$ centered on
$p$. It follows that $\ray{\sigma}{\xi}$ and $\ray{\tau}{\eta}$ have
infinite Hausdorff distance in $\T$.

The second sentence follows from the first, by dividing each line into two rays.
\end{proof}

\paragraph{Remark} Combining results of Masur mentioned above, one can show that even more is
true: two cobounded geodesic rays which have finite Hausdorff distance are asymptotic, meaning
that as they go to $\infinity$, the distance between the rays approaches zero. To see why, as
mentioned earlier Masur proves that if $\ray{\sigma}{\eta}$ is cobounded then $\eta$ is uniquely
ergodic. Furthermore, two rays $\ray{\sigma}{\eta}$, $\ray{\tau}{\eta}$ with uniquely
ergodic endpoint $\eta$ are asymptotic, according to \cite{Masur:UniquelyErgodic}.

\subsection{Singular {\small\rm SOLV} spaces}
\label{SectionSolv}

Consider a geodesic $g = \geodesic{\xi}{\eta}$ in $\T$, and let $\S_g \to g$ be the canonical
marked Riemann surface bundle over $g$, obtained by pulling back the canonical marked
Riemann surface bundle $\S \to \T$. Topologically we identify $\S_g = S \cross g$. Choosing a
transverse pair of measured foliations $\F_x, \F_y$ representing $\xi,\eta$ respectively, we
have $g(t) = \sigma(e^{-t} \F_x, e^{t} \F_y)$. Let $\abs{dy}$ be the transverse measure on
the horizontal measured foliation $\F_x$ and let $\abs{dx}$ be the transverse measure on the
vertical measured foliation $\F_y$. We may assume that the pair $\F_x,\F_y$ is normalized,
meaning that the Euclidean area equals~$1$:
$$\norm{q(\F_x,\F_y)} = \int_S \abs{dx} \cross \abs{dy} = 1
$$
and hence for all $t \in \reals$ the pair $e^{-t} \F_x, e^t \F_y$ is normalized:
$$\norm{q(e^{-t}\F_x,e^{t}\F_y)} = \int_S \abs{e^t dx} \cross \abs{e^{-t} dy} = 1
$$
Note that the singular Euclidean metric on each fiber $\S_{g(t)}$, may be expressed as
$$ds_\sigma^2 = e^{2t} \abs{dx}^2 + e^{-2t} \abs{dy}^2
$$
Define the \emph{singular \solv\ metric} on $\S_g$ to be the singular Riemannian metric
given by the formula:
$$ds_g^2 = e^{2t} \abs{dx}^2 + e^{-2t} \abs{dy}^2 + dt^2
$$
We use the notation $\S^\solv_g$ to denote $\S_g$ equipped with this metric. The
universal cover of $\S_g$ is the canonical \Poincare\ disc bundle $\H_g$ over $g$, and
lifting the singular \solv\ metric from $\S^\solv_g$ to $\H_g$ we obtain a singular
\solv\ space denoted $\H^\solv_g$. The singular locus of $\S^\solv_g = S \cross
g$ is the union of the singular lines $s \cross g$, one for each singularity $s$ of
the pair $\F_x, \F_y$. Away from the singular lines, $\S^\solv_g$ and $\H^\solv_g$
are locally modelled on 3--dimensional \solv--geometry. On
each singular line the metric is locally modelled by gluing together several copies of the
half-plane $y\ge 0$ in $\solv$--geometry.

\subsection{Comparing hyperbolic and singular Euclidean structures}
\label{SectionComparisons}

Given $\sigma \in \T$, the Riemann surface $\S_\sigma$ has several important metrics in its
conformal class: a unique hyperbolic metric; and one singular Euclidean metric for each $q\in
\QD_\sigma$. These lift to the universal cover $\H_\sigma$. Given
$\sigma,\tau \in \T$, if each Riemann surface $\S_\sigma$, $\S_\tau$ is
given either its unique hyperbolic metric or one of its singular
Euclidean metrics, then for any marked map $\phi \from\S_\sigma \to
\S_\tau$, each lift $\wt\phi\from \H_\sigma \to \H_\tau$ is a quasi-isometry. We are interested in
how the quasi-isometry constants of $\wt\phi$ compare to the \Teichmuller\ distance
$d(\sigma,\tau)$, although we need only the crudest estimates. Proposition~\ref{PropTeichQuasi}
shows how to bound the quasi-isometry constants in terms of $d(\sigma,\tau)$. Part~1 of this
proposition was first proved by Minsky in \cite{Minsky:endinglaminations}, Lemma 3.3; we give a
quicker proof using Lemma~\ref{LemmaFamilyOfMetrics}. 

\begin{proposition} 
\label{PropTeichQuasi} 
For each bounded subset $\B \subset \Mod$ and each $r>0$ there exists $K\ge 1, C  \ge 0, A \ge
0$ such that the following hold:
\begin{enumerate}
\item Suppose that $\sigma,\tau \in \T$ are each $\B$--cobounded and $d(\sigma,\tau) \le r$. Let
$f_{\sigma \tau} \from \S_\sigma\to \S_\tau$ be the canonical marked map
$\S_\sigma=S\cross\sigma \to S\cross\tau=\S_\tau$. If we impose on $\S_\sigma$ and $\S_\tau$
either the hyperbolic metric or the singular Euclidean metric associated to some normalized
quadratic differential, then any lift $\wt f_{\sigma \tau} \from
\H_\sigma \to \H_\tau$ of $f_{\sigma\tau}$ is a $K,C$ quasi-isometry.
\item Let $\sigma_i \in \T$, $i=1,2,3$, be $\B$--cobounded and have pairwise distances $\le r$,
let metrics be imposed on $\S_{\sigma_i}$ as above, and let $f_{ij} \from \S_{\sigma_i}
\to\S_{\sigma_j}$, etc.\ be the marked maps as above, with $K,C$--quasi-isometric lifts $\wt f_{ij}
\from\H_{\sigma_i} \to \H_{\sigma_j}$. If $\wt f_{13}$ is the unique lift of $f_{13}$ such
that
$$\bdy \wt f_{23} \composed \bdy \wt f_{12} = \bdy \wt f_{13},
$$ 
then
$$\dsup(\wt f_{23} \composed \wt f_{12}, \wt f_{13}) 
\le A.
$$
\end{enumerate}
\end{proposition}

\begin{proof} Part (1) is an easy consequence of Lemma~\ref{LemmaFamilyOfMetrics}, as
follows. Choose a compact subset $\A \subset\T$ whose image in $\Mod$ covers $\B$ and such
that over any point of $\B$ there exists a point $\sigma \in \A$ such that $\B_\T(\sigma,r)
\subset\A$. It follows that the points $\sigma,\tau$ in (1) may be translated to lie in $\A$.
Identifying $\S_\A$ diffeomorphically with $S \cross \A$, compactness of $\A$ produces a
compact family of hyperbolic metrics on $S$, and compactness of 
the restriction of $\QD^1$ to $\A$ produces a
compact family of singular Euclidean metrics. Now apply Lemma~\ref{LemmaFamilyOfMetrics}.

For part (2), note that by compactness of $\A$ and of the compactness of
the restriction of
$\QD^1$ to $\A$, there exists a
uniform $\delta$ such that any hyperbolic metric and any normalized singular Euclidean
structure determined by an element $\sigma\in \A$ has a $\delta$--hyperbolic universal cover.
Part (2) is now a direct consequence of Lemma~\ref{LemmaBdyValueUniqueness}.
\end{proof}

\section{Convex cocompact subgroups of $\Isom(\T)$}
\label{SectionConvexCocompactSubgroups}

\subsection{Variations of convex cocompactness}
\label{SectionVariations}

Given a proper, geodesic metric space $X$, a subset $L \subset X$ is \emph{quasiconvex} if
there exists $A \ge 0$ such that every geodesic segment in $X$ with endpoints in $L$ is
contained in the $A$--neighborhood of $L$.

When $G$ is a finitely generated, discrete subgroup of the isometry group of $\hyp^n$, it is
well known that the following properties of $G$ are all equivalent to each other:
\begin{description} 
\item[Orbit Quasiconvexity] Any orbit of $G$ is a quasiconvex subset of
$\hyp^n$.
\item[Single orbit quasiconvexity] There exists an orbit of $G$ which is quasiconvex in
$\hyp^n$.
\item[Convex cocompact] $G$ acts cocompactly on the convex hull of its limit set
$\Lambda$. 
\end{description}
Moreover, these properties imply that $G$ is word hyperbolic, and there
is a continuous $G$--equivariant embedding of the Gromov boundary $\bdy
G$ into $\bdy\hyp^n$ whose image is the limit set $\Lambda$. Similar
facts hold for finitely generated groups acting discretely on any Gromov
hyperbolic space, for example finitely generated subgroups of Gromov
hyperbolic groups.

In this section we prove Theorem~\ref{TheoremQuasiconvex}, which is a list of similar
equivalences for finitely generated subgroups of the isometry group of the \Teichmuller\
space $\T$ of~$S$. In this case the entire isometry group $\Isom(\T)$ acts discretely on
$\T$, and in fact by Royden's Theorem \cite{Royden:IsomTeichReport}, \cite{Ivanov:Automorphisms}
the canonical homomorphism $\MCG \to \Isom(\T)$ is an isomorphism, except in genus~2 where the
kernel is cyclic of order~2.

Although $\T$ fails to be negatively curved in any reasonable sense, nevertheless one can say
that it behaves in a negatively curved manner as long as one focusses only on cobounded
aspects. This, at least, is one way to interpret the projection properties introduced by
Minsky in \cite{Minsky:quasiprojections} and further developed by Masur and Minsky in
\cite{MasurMinsky:complex1}. Given a $\B$--cobounded geodesic $g$ in $\T$, Minsky's
projection property says that a closest point projection map of $\T$ onto $g$ behaves in
a negatively curved manner, such that the quality of the negative curvature depends only on
$\B$. See Theorem~\ref{TheoremMinskyContraction} for the precise statement.

For a finitely generated subgroup $G\subset \Isom(\T)$ we can obtain equivalences as above,
as long as we tack on an appropriate uniform coboundedness property; in some cases the desired
property comes for free by uniform coboundedness of the action of $G$ on any of its orbits.

First we have some properties of $G$ which are variations on orbit quasiconvexity:
\begin{description}
\item[Orbit quasiconvexity] Any orbit of $G$ is quasiconvex in $\T$.
\item[Single orbit quasiconvexity] There exists an orbit of $G$ 
that is quasiconvex in
$\T$.
\item[Weak orbit quasiconvexity] 
There exists a constant $A$ and an orbit $\O$ of $G$, and for each $x,y
\in \O$ there exists a geodesic segment $[x',y']$ in $\T$, such that
$d(x,x') \le A$, $d(y,y') \le A$, and $[x',y']$ is in the
$A$--neighborhood of $\O$.
\end{description}
The latter is a more technical version of orbit quasiconvexity which is quite useful in
several settings. 

Another property of $G$ is a version of convex cocompactness, into which we
incorporate the hyperbolicity properties mentioned above:
\begin{description}
\item[Convex cocompact] The group $G$ is word hyperbolic, and there exists a continuous
$G$--equivariant embedding $f_\infinity \from \bdy G \to \PMF$ with image $\Lambda_G$, such that
$\Lambda_G \cross \Lambda_G - \Delta \subset \P\FP$, and the following holds. Letting 
$$\WHull_G = \union \{\geodesic{\zeta}{\zeta'} \suchthat \zeta\ne\zeta' \in\Lambda_G\}
$$
be the \emph{weak hull} of $\Lambda_G$, if $f \from G \to \WHull_G$ is any
$G$--equivariant map, then $f$ is a quasi-isometry and the map $\overline f = f \union
f_\infinity \from G
\union
\bdy G \to \WHull_G \union \Lambda_G$ is continuous. 
\end{description}

In this definition, $\WHull_G$ is metrized by restricting the
\Teichmuller\ metric on~$\T$, which \emph{a posteriori} has the effect
of making $\WHull_G$ into a quasigeodesic metric space. The definition
implies that $G$ acts cocompactly on $\WHull_G$: since $\Lambda_G
\cross \Lambda_G - \Delta$ is a closed subset of $\P\FP$ it follows that $\WHull_G$ is a closed
subset of $\T$; and since $G$ acts coboundedly on itself it follows that $G$ acts coboundedly on
$\WHull_G$; thus, the image of $\WHull_G$ in moduli space is closed and bounded, hence compact.

\subsection{Properties of convex cocompact subgroups}
\label{SectionProperties}

In this section we prove several properties of convex cocompact subgroups of $\Isom(\T)$ which
are analogues of well known properties in $\Isom(\hyp^n)$.

\begin{proposition} 
\label{PropCnvxCbddPsA}
Every infinite order element $g$ of a convex cocompact subgroup $G\subgroup \Isom(\T) \approx
\MCG$ is a pseudo-Anosov mapping class. 
\end{proposition}

\begin{proof} Any infinite order element of a word hyperbolic group has 
source--sink dynamics on
its Gromov boundary, and so $g$ has source--sink dynamics on $\bdy G\approx \Lambda_G$. It follows
that $g$ has an axis in $\WHull_G$. But the elements of $\Isom(\T) \approx \MCG$ having an axis in
$\T$ are precisely the pseudo-Anosovs \cite{Bers:ThurstonTheorem}.
\end{proof}

The following is a consequence of work of McCarthy and Papadoupolos \cite{McCarthyPapa:dynamics}.

\begin{proposition}
\label{PropCnvxCbddLimitSet}
If $G$ is a convex cocompact subgroup of $\Isom(\T)$ then:
\begin{enumerate}
\item $\Lambda_G$ is the smallest nontrivial closed subset of $\overline\T = \T \union \PMF$
invariant under~$G$.
\item The action of $G$ on $\PMF\setminus\Lambda_G$ is properly discontinuous.
\end{enumerate}
\end{proposition}

\begin{proof} The Gromov boundary of a word hyperbolic group is the closure of the fixed points
of infinite order elements in the group, and so by Proposition~\ref{PropCnvxCbddPsA} the set
$\Lambda_G$ is the closure of the fixed points of the pseudo-Anosov elements of $G$. Item (1) now
follows from Theorem 4.1 of \cite{McCarthyPapa:dynamics}. 

To prove (2), let
$$Z(\Lambda) = \{\zeta \in \PMF \suchthat \text{there exists}\quad \zeta' \in \Lambda
\quad\text{such that}\quad i(\zeta,\zeta') = 0 \}
$$
Theorem 6.16 of \cite{McCarthyPapa:dynamics} says that $G$ acts properly discontinuously on
$\PMF-Z(\Lambda)$, and so it suffices to prove that $\Lambda=Z(\Lambda)$. Each point
$\zeta'\in\Lambda$ is the ideal endpoint of a cobounded geodesic ray, which implies that
$\zeta'$ is uniquely ergodic and fills the surface \cite{Masur:exchange}, and so if
$i(\zeta,\zeta')=0$ then $\zeta=\zeta'$.
\end{proof}

\paragraph{Remark} One theme of \cite{McCarthyPapa:dynamics} is that for a general finitely
generated subgroup $G \subgroup \MCG$, there are several different types of ``limit sets'' for
the action of $G$ on $\PMF$. Assuming that $G$ contains a pseudo-Anosov element, the two sets
mentioned in the proof above play key roles in \cite{McCarthyPapa:dynamics}: $\Lambda(G)$ which
is the closure of the fixed points of pseudo-Anosov elements of the subgroup, and is also the
smallest nontrivial closed $G$--invariant subset; and the set $Z(\Lambda(G))$. What we have proved
is that for a convex cocompact subgroup $G$, these two sets are identical. Henceforth we refer to
$\Lambda_G$ as \emph{the limit set} for the action of $G$ on $\PMF$.

\medskip

The analogue of the following result is true for convex cocompact discrete subgroups of
$\hyp^n$, as well as for word hyperbolic groups \cite{Arzhantseva:quasiconvex}; the proof here is
similar.

\begin{proposition} 
\label{PropConvCbddComm}
Let $G$ be a convex cocompact subgroup of $\Isom(\T)$, and let $N_G$ and
$\Comm_G$ be the normalizer and the relative commensurator of $G$ in
$\Isom(\T)$. Then each of the inclusions $G \subgroup N_G \subgroup
\Comm_G$ is of finite index, and we have $\Comm_G = \Stab(\Lambda_G) =
\Stab(\WHull_G)$.
\end{proposition}

\begin{proof} Let $\Lambda_G$ be the limit set of $G$, with weak hull
$\WHull_G$, and note that we trivially have $\Stab(\WHull_G) =
\Stab(\Lambda_G)$. 

Note that $\Stab(\WHull_G)$ acts properly on $\WHull_G$. Indeed,
$\Isom(\T)$ acts properly on $\T$, and so any subgroup of $\Isom(\T)$
acts properly on any subset of $\T$ which is invariant under that
subgroup. Since $G\subset \Stab(\WHull_G)$, and since $G$ acts
cocompactly on $\WHull_G$, it follows that $G$ is contained with finite
index in $\Stab(\WHull_G)$. This implies that $\Stab(\WHull_G) \subset
\Comm_G$.  To complete the proof we only have to
prove the reverse inclusion $\Comm_G \subset
\Stab(\WHull_G)$. 

Given $g \in \Isom(\T)$, suppose that $g \in \Comm_G$, and choose finite
index subgroups $H,K \subgroup G$ such that $g^\inv H g = K$. By the
definition of convex cocompactness it follows that $\WHull_H =
\WHull_G = \WHull_K$. Since $g(\WHull_K) = \WHull_H$ it follows that $g \in \Stab(\WHull_G)$.
\end{proof}

\paragraph{Remark}
Another natural property for subgroups $G \subgroup \MCG$ is
quasiconvexity with respect to the word metric on $\MCG$. It seems
possible to us that this is not equivalent to orbit quasiconvexity of
$G$ in $\Isom(\T)$. Masur and Minsky \cite{MasurMinsky:unstable} give an
example of an infinite cyclic subgroup of $\Isom(\T)$ which is not orbit
quasiconvex, and yet this subgroup is quasi-isometrically embedded in
$\MCG$ \cite{FarbLubotzkyMinsky}; it may also be quasiconvex in $\MCG$,
but we have not investigated this.

\subsection{Equivalence of definitions: Proof of 
Theorem~\ref{TheoremQuasiconvex}}
\label{SectionEquivalence}

Here is our main result equating the various quasiconvexity properties
with convex cocompactness:

\ppar
\noindent
{\bf Theorem \ref{TheoremQuasiconvex}}\qua
{\sl If $G$ is a finitely generated subgroup of $\Isom(\T)$, the following are
equivalent:
\begin{enumerate}
\item Orbit quasiconvexity
\item Single orbit quasiconvexity
\item Weak orbit quasiconvexity
\item Convex cocompactness
\end{enumerate}
}

\ppar
Because of this theorem we are free to refer to ``quasiconvexity'' or ``convex
cocompactness'' of $G$ without any ambiguity.

\begin{proof}[Proof of Theorem \ref{TheoremQuasiconvex}]
 The key ingredients in the proof are results of Minsky from
\cite{Minsky:quasiprojections} concerning projections from balls and horoballs in $\T$ to
geodesics in $\T$, and results of Masur--Minsky
\cite{MasurMinsky:complex1} characterizing $\delta$--hyperbolicity of
proper geodesic metric spaces in terms of projections properties to
paths.

To begin with, note that the implications $(1)\Rightarrow (2)\Rightarrow 
(3)$ are obvious.  We now prove that $(3)\Rightarrow (1)$. 

Suppose we have an orbit $\O$ of $G$ and a constant $A$, and for each
$x,y \in \O$ we have two points $x',y' \in \T$, endpoints of a unique
geodesic segment $[x',y']$ in $\T$, such that $d(x,x') \le A$, $d(y,y')
\le A$, and $[x',y'] \subset N_A(\O)$. The set $\O$ maps to a single
point in $\Mod$ and so the projection of $N_A(\O)$ to $\Mod$ is a
bounded set $\B$. It follows that each $[x',y']$ is $\B$--cobounded. Now
consider an arbitrary orbit $\O_1$ of $G$; we must prove that $\O_1$ is
quasiconvex in $\T$. The orbits $\O, \O_1$ have finite Hausdorff
distance $C$ in $\T$. Given $x_1,y_1 \in\O_1$, choose $x,y \in \O$
within distance $C$ of $x_1,y_1$, respectively, and consider the
geodesic segment $[x',y']$ and the piecewise geodesic path $$\gamma =
[x',x] * [x,x_1] * [x_1,y_1] * [y_1,y] * [y,y'] $$ Of the five
subsegments of $\gamma$, all but the middle subsegment have length $\le
\Max\{A,C\}$, and it follows that $\gamma$ is a $(1,D)$--quasigeodesic in $\T$, with $D$
depending only on $A,C$. Since the geodesic $[x',y']$ is $\B$--cobounded we can apply the
following result of Minsky \cite{Minsky:quasiprojections} to obtain $\delta$,
depending only on $\B$ and $D$, such that $\gamma \subset N_\delta [x',y']$. 

\begin{theorem}[Stability of cobounded geodesics]
\label{TheoremMinskyStability}
For any bounded subset $\B$ of $\Mod$ and any $K \ge 1, C \ge 0$ there exists $\delta \ge 0$
such that if $\gamma$ is a $K,C$ quasigeodesic in $\T$ with endpoints $x,y$, and if $[x,y]$ is
$\B$--cobounded, then $\gamma \subset N_\delta[x,y]$.
\qed\end{theorem}

It follows that $[x_1,y_1] \subset \gamma \subset N_{\delta+A}\O \subset N_{\delta+A+C} \O_1$,
proving quasiconvexity of $\O_1$ in $\T$.

\paragraph{Weak orbit quasiconvexity implies convex cocompactness} 
Fix an orbit $\O$ of $G$, and so $\O$ is quasiconvex in $\T$. Let $\G$ be the set of all geodesic
segments, rays, and lines that are obtained as pointwise limits of sequences of geodesics with
endpoints in $\O$. Let $\union\G \subset \T$ be the union of the elements of $\G$. The left
action of $G$ on $\O$ is evidently cobounded. By quasiconvexity of $\O$ it follows that the
action of $G$ on the union of geodesic segments with endpoints in $\O$ is cobounded, which
implies in turn that the action of $G$ on $\union\G$ is cobounded. Since $\union\G$ is closed and
$\T$ is locally compact, it follows that the $G$ action on $\union\G$ is cocompact. The set
$\union\G$ therefore projects to a compact subset of $\Mod$ which we denote $\B$. All geodesics
in $\G$ are therefore $\B$--cobounded.

Let $\union\G$ be equipped with the restriction of the \Teichmuller\ metric. Note that while
$\union\G$ is not a geodesic metric space, it is a quasigeodesic metric space: there exists $A
\ge 0$ such that any $x,y \in \union\G$ are within distance $A$ of points $x',y' \in
\O \subset \union\G$, and the geodesic $[x',y']$ is contained in $\union\G$. 

To prepare for the proof that $G$ is word hyperbolic, fix a finite generating set for $G$ with
Cayley graph $\Gamma$, and fix a $G$--equivariant map $f \from \Gamma \to \union\G$ taking the
vertices of $\Gamma$ to $\O$ and taking each edge of $\Gamma$ to an element of
$\G$. Since $G$ acts properly and coboundedly on both $\Gamma$ and $\union\G$, and since both are
quasigeodesic metric spaces, it follows that the equivariant map $f$ is a quasi-isometry between
$\Gamma$ and $\union\G$; pick a coarse inverse $F \from\union\G \to \Gamma$.

By definition the group $G$ is word hyperbolic if and only if the Cayley graph $\Gamma$ is
$\delta$--hyperbolic for some $\delta \ge 0$. Our proof that $G$ is word hyperbolic will use a
result of Masur and Minsky, Theorem 2.3 of \cite{MasurMinsky:complex1}:

\begin{theorem}
\label{TheoremMMHyperbolic}
Let $X$ be a geodesic metric space and suppose that there is a set of paths $\P$ in
$X$ with the following properties:
\begin{description}
\item[Coarse transitivity] 
There exists $C \ge 0$ such that for any $x,y \in X$ with $d(x,y)
\ge C$ there is a path in $\P$ joining $x$ and $y$.
\item[Contracting projections:] There exist $a,b,c > 0$, and for each path $\gamma \from I \to
X$ in $\P$ there exists a map $\pi \from X \to I$ such that:
\begin{description}
\item[Coarse projection] For each $t \in I$ we have $\diam\left(\gamma[t,\pi(\gamma t)]\right)
\le c$.
\item[Coarse lipschitz] If $d(x,y) \le 1$ then $\diam\left(\gamma[\pi x, \pi y]\right) \le c$.
\item[Contraction] If $d(x,\gamma(\pi x)) \ge a$ and $d(x,y) \le b \cdot d(x,\gamma(\pi x))$
then $$\diam\left(\gamma[\pi x, \pi y] \right) \le c $$
\end{description}
\end{description}
Then $X$ is $\delta$--hyperbolic for some $\delta \ge 0$.
\end{theorem}

To prove that $G$ is $\delta$--hyperbolic we take $\P$ to be the set of geodesic segments
in $G$, and we look at the set of paths $f \composed \P = \{f\composed\gamma \suchthat
\gamma\in \P\}$ in $\union\G$. Using some results of Minsky \cite{Minsky:quasiprojections},
we will show that $f \composed \P$ satisfies the hypotheses of Theorem
\ref{TheoremMMHyperbolic}. Then we shall pull the hypotheses back to $\P$ and apply Theorem
\ref{TheoremMMHyperbolic}.

The first result of Minsky that we need is the main theorem of \cite{Minsky:quasiprojections}:

\begin{theorem}[Contraction Theorem]
\label{TheoremMinskyContraction} 
For every bounded subset $\B$ of $\Mod$ there exists $c>0$ such that if $\gamma$ is any
$\B$--cobounded geodesic in $\T$ then the closest point projection $\T \to \gamma$ satisfies
the $(a,b,c)$ contracting projection property with $(a,b)=(0,1)$.
\end{theorem}

In our context, where we have a uniform $\B$ such that each geodesic in $\G$ is
$\B$--cobounded, it follows that there is a uniform $c$ such that each geodesic in $\G$
satisfies the $(0,1,c)$ contracting projection property.

Now consider $\gamma = [x_0,x_1,\ldots,x_n]$ a geodesic in the Cayley graph $\Gamma$, mapping
via $f$ to a piecewise geodesic $f\gamma = [f x_0,f x_1] \union \cdots \union [f x_{n-1},f
x_n]$ in $\union\G$, with each subsegment $[f x_i,f x_{i+1}]$ an element of $\G$. It follows
that $f\gamma$ is a $K,C$ quasigeodesic in $\T$, for $K \ge 1, C \ge 0$ independent of the
given geodesic in $\Gamma$. The $\T$--geodesic $[f x_0,f x_n]$ is $\B$--cobounded. Applying
Theorem \ref{TheoremMinskyStability} it follows that $f\gamma \subset N_D[f x_0,f x_n]$, where
$D$ depends only on $\B,K,C$. As noted above, closest point projection from $\T$ onto $[f
x_0, f x_n]$ satisfies the $(0,1,c)$ contracting projection property. From this it follows
that closest point projection $\pi \from \T \to f\gamma$ satisfies the $(a',b',c')$
contraction property where $(a',b',c')$ depend only on $\B,K,C$. Now define the projection
$\Gamma \to \gamma$ to be the composition $\Gamma\xrightarrow{f} \union\G \xrightarrow{\pi}
f\gamma \xrightarrow{F}\Gamma \to \gamma$ where the last map is closest point
projection in $\Gamma$. This composition clearly satisfies the $(a'',b'',c'')$ projection
property where $(a'',b'',c'')$ depend only on $(a',b',c')$ and the quasi-isometry constants
and coarse inverse constants for $f,F$.

Geodesics in $\Gamma$ are clearly coarsely transitive, and applying Theorem
\ref{TheoremMMHyperbolic} it follows that $G$ is word hyperbolic. This means that geodesic
triangles in $\Gamma$ are uniformly thin, and it implies that for each $K,C$ there is a
$\delta$ such that $K,C$ quasigeodesic triangles in $\Gamma$ are $\delta$--thin. Applying the
quasi-isometry between $\Gamma$ and $\union\G$, it follows that there is a uniform $\delta$
such that for each $x,y,z \in \O$ the geodesic triangle $\triangle[x,y,z]$ in $\union\G$ is
$\delta$--thin; we fix this $\delta$ for the arguments below.

Now we turn to a description of the ``limit set'' $\Lambda \subset \PMF$ of $G$, with the
ultimate goal of identifying it with the Gromov boundary $\bdy G$. 

Each geodesic ray in $\G$ has the form $\ray{x}{\eta}$, for some $x \in \O$, $\eta \in\PMF$;
define $\Lambda \subset \PMF$ be the set of all such points $\eta$, over all geodesic rays in
$\G$. The set $\Lambda$ is evidently $G$--equivariant. 

\subparagraph{Fact 1} For any $x \in \O$, $\eta \in \Lambda$, the ray $\ray{x}{\eta}$ in
$\T$ is an element of~$\G$.

\medskip
To prove this, by definition of $\Lambda$ there exists a ray
$\ray{y}{\eta}$ in $\G$ for some $y \in \O$. Choose a sequence
$y_1,y_2,\ldots \in \O$ staying uniformly close to $\ray{y}{\eta}$ and
going to infinity. Pass to a subsequence so that the sequence of
segments $[x,y_n]$ converges to some ray $\ray{x}{\eta'} \in \G$; it
suffices to show that $\eta'=\eta$. Since $x$ is fixed and the points
$y_n$ stay uniformly close to $\ray{y}{\eta}$, it follows by Theorem
\ref{TheoremMinskyStability} that the segments $[x,y_n]$ stay uniformly
close to $\ray{y}{\eta}$, and so $\ray{x}{\eta'}$ is in a finite
neighborhood of $\ray{y}{\eta}$. The reverse inclusion, that
$\ray{y}{\eta}$ is in a finite neighborhood of $\ray{x}{\eta'}$, is a
standard argument: as points move to infinity in $\ray{x}{\eta'}$ taking
bounded steps, uniformly nearby points move to infinity in
$\ray{y}{\eta}$ also taking bounded steps, and thus must come uniformly
close to an arbitrary point of $\ray{y}{\eta}$. This shows that the rays
$\ray{x}{\eta'}$, $\ray{y}{\eta}$ have finite Hausdorff distance, and
applying Lemma \ref{LemmaEndUniqueness} (End Uniqueness) shows that
$\eta=\eta'$.

Note that in the proof of Fact 1 we have established a little more, namely that for any $x,y
\in \O$ and $\eta \in \Lambda$ the rays $\ray{x}{\eta}$ and $\ray{y}{\eta}$ have finite Hausdorff
distance. This will be useful below.

\subparagraph{Fact 2} For any $\eta \ne \zeta \in \Lambda$ there exists a line
$\geodesic{\eta}{\zeta}$ contained in $\G$.

\medskip
From Fact 2 it immediately follows that $\Lambda \cross \Lambda - \Delta
\subset \P\FP$, that the weak hull $\WHull_G$ of $\Lambda$ is defined,
and that $G$ acts coboundedly on $\WHull_G$, since $G$ acts coboundedly
on $\union\G$.

To prove Fact 2, pick a point $x \in \O$, and note that by Fact 1 we have two rays
$\ray{x}{\eta}$, $\ray{x}{\zeta}$ in $\G$. Pick a sequence $y_n \in \O$ staying uniformly close to
$\ray{x}{\eta}$ and going to infinity, and a sequence $z_n \in \O$ staying uniformly close to
$\ray{x}{\zeta}$ and going to infinity. We have a sequence of triangles $[x,y_n,z_n]$ in $\G$, all
$\delta$--thin. Applying Theorem \ref{TheoremMinskyStability} there is a $D$ such that the sides
$[x,y_n]$ are contained in the $D$--neighborhood of $\ray{x}{\eta}$, and the sides $[x,z_n]$ are
contained in the $D$--neighborhood of $\ray{x}{\zeta}$. Each side $[y_n,z_n]$, being contained in
the $\delta$--neighborhood of $[x,y_n] \union [x,z_n]$, is therefore contained in the
$D+\delta$--neighborhood of $\ray{x}{\eta} \union \ray{x}{\zeta}$. 

We claim that the point $x$ is uniformly close to the segments $[y_n,z_n]$. If not, then from
uniform thinness of the triangles $[x,y_n,z_n]$ it follows that there are points $y'_n \in
[x,y_n]$ and $z'_n \in [x,z_n]$ such that the segments $[x,y'_n]$ and $[x,z'_n]$ get arbitrarily
long while the Hausdorff distance between them stays uniformly bounded.
This implies that there are sequences $y''_n \in\ray{x}{\eta}$ going to infinity and $z''_n \in
\ray{x}{\zeta}$ going to infinity such that the Hausdorff distance between the segments
$[x,y''_n]$ and $[x,z''_n]$ stays uniformly bounded, which implies in turn that the rays
$\ray{x}{\eta}$ and $\ray{x}{\zeta}$ have finite Hausdorff distance. Applying End
Uniqueness~\ref{LemmaEndUniqueness}, it follows that $\eta=\zeta$, contradicting the hypothesis
of Fact~2, and the claim follows.

Passing to a subsequence and applying Ascoli--Arzela it follows that $[y_n,z_n]$ converges to a
line in $\G$. One ray of this line is Hausdorff close to $\ray{x}{\eta}$ and so has endpoint
$\eta$, and the other ray is Hausdorff close to $\ray{x}{\zeta}$ and so has endpoint $\zeta$, by
End Uniqueness. We therefore have $\lim[y_n,z_n] = \geodesic{\eta}{\zeta}$, completing the proof
of Fact~2.

\medskip

Now we define a map $f_\infinity \from \bdy G \to \Lambda$. Recall that the relation of
finite Hausdorff distance is an equivalence relation on geodesic rays in the Cayley
graph~$\Gamma$ of $G$, and $\bdy G$ is the set of equivalence classes. Consider then a point
$\xi\in \bdy G$ represented by two geodesic rays $[x_0,x_1,\ldots)$ and $[y_0,y_1,\ldots)$ with
finite Hausdorff distance in $\Gamma$. These map to piecewise geodesic, quasigeodesic rays $\rho
= [fx_0,fx_1] \union [fx_1,fx_2] \union \cdots$ and $\sigma = [fy_0,fy_1]\union [fy_1,fy_2]
\union\cdots$ with finite Hausdoff distance in $\union\G$. The sequence of geodesic segments
$[fx_0,fx_n]$ in $\G$ has a subsequence converging to some ray $\ray{fx_0}{\zeta}$ in $\G$, and
$[fy_0,fy_n]$ has a subsequence converging to some ray $\ray{fy_0}{\zeta'}$ in $\G$. To obtain a
well defined map $\bdy G \to \Lambda$ it suffices to prove that $\zeta=\zeta'$, and then we can
set $f_\infinity(\xi) = \zeta$. 

To prove that $\zeta=\zeta'$ it suffices, by End Uniqueness \ref{LemmaEndUniqueness}, to prove
that the rays $\ray{fx_0}{\zeta}$ and $\ray{fy_0}{\zeta'}$ have finite Hausdorff distance in
$\T$. Since the piecewise geodesic rays $\rho,\sigma$ have finite Hausdorff distance in $\T$, it
suffices to prove that $\rho$ has finite Hausdorff distance from $\ray{fx_0}{\zeta}$, and
similarly $\sigma$ has finite Hausdorff distance from $\ray{fy_0}{\zeta'}$. Consider a point $p
\in\rho$. For sufficiently large $n$ we have $p \in \rho_n = [fx_0,fx_1] \union\cdots\union
[fx_{n-1},fx_n]$. Applying Theorem \ref{TheoremMinskyStability} there is a uniform constant $D$
such that $\rho_n \subset N_D([fx_0,fx_n])$, and so $p$ is within distance $D$ of some point in
$[fx_0,fx_n]$. Since $\ray{fx_0}{\zeta}$ is the pointwise limit of $[fx_0,fx_n]$ as
$n\to\infinity$ it follows that $p$ is within a uniformly bounded distance of
$\ray{fx_0}{\zeta}$. This shows that $\rho$ is within a finite neighborhood of
$\ray{fx_0}{\zeta}$. The reverse inclusion is a standard argument: as points move along $\rho$
towards the end taking bounded steps, uniformly nearby points move along $\ray{fx_0}{\zeta}$
towards the end also taking bounded steps, and thus must come uniformly close to some point of
$\ray{fx_0}{\zeta}$.

Hence $f_\infinity \from \bdy G \to \Lambda$ is well defined. Observe that a similar argument
proves a little more: if $x_i \in G$ converges to $\xi \in \bdy G$ then the segments
$[fx_0,fx_i]$ converge in the compact--open topology to the ray $\ray{fx_0}{f\xi}$; details are
left to the reader.

We now turn to verifying required properties of $f_\infinity$.

To see that $f_\infinity$ is surjective, consider a point $\eta \in\Lambda$ and pick a ray
$\ray{x}{\eta}$ in $\G$. It follows that $\rho = F\left(\ray{x}{\eta}\right)$ is a
quasigeodesic ray in $\Gamma$. Since $\Gamma$ is $\delta$--hyperbolic it follows that $\rho$ has
finite Hausdorff distance from some geodesic ray $\rho'$ in $\Gamma$, with endpoint $\zeta'
\in\bdy G$. As shown above, $f(\rho')$ has finite Hausdorff distance from some geodesic ray
$\ray{x'}{f_\infinity \zeta'}$. Since $f, F$ are coarse inverses it follows that
$\ray{x}{\eta}$ has finite Hausdorff distance from $\ray{x'}{f_\infinity \zeta'}$, and so by End
Uniqueness it follows that $\eta = f_\infinity\zeta'$.

To see that $f_\infinity$ is injective, consider two points $\eta,\zeta\in \bdy G$ and suppose
that $f_\infinity(\eta)=f_\infinity(\zeta)$; let $\xi \in \Lambda$ be this point. Pick rays
$\rho,\sigma$ in $\Gamma$ representing $\eta,\zeta$ respectively. As we have just seen, the images
$f(\rho)$, $f(\sigma)$ have finite Hausdorff distance in $\T$ to rays $\ray{y}{\xi}$,
$\ray{z}{\xi}$ in $\G$, respectively. As noted at the end of the proof of Fact 1, the rays
$\ray{y}{\xi}$ and $\ray{z}{\xi}$ have finite Hausdorff distance in $\T$; applying the coarse
inverse $F$ it follows that $\rho,\sigma$ have finite Hausdorff distance in $\Gamma$ and
therefore $\eta=\zeta$.

We have shown that $f_\infinity$ is a bijection between $\bdy G$ and $\Lambda$. We want to prove
that $f_\infinity$ is a homeomorphism, and that the extension $\overline f = f \union f_\infinity
\from G \union \bdy G \to \overline\T = \T \union \PMF$ is continuous. For this purpose first we
establish:

\subparagraph{Fact 3} $\Lambda$ is a closed subset of $\PMF$, and therefore compact.

\medskip
To prove this, choose a sequence $\zeta_n \in \Lambda$ so that $\lim\zeta_n =
\zeta_\infinity$ in $\PMF$; we must prove that $\zeta_\infinity \in \Lambda$. Choose a point
$x \in \O$, and apply Fact 1 to obtain rays $\ray{x}{\zeta_n}$. Passing to a subsequence these
converge to a limiting ray $\lim\ray{x}{\zeta_n}=\ray{x}{\zeta'_\infinity}$ in $\G$, and so
$\zeta'_\infinity \in \Lambda$. Looking in the unit tangent bundle of $\T$ at the point $x$ it
follows that $\lim\zeta_n = \zeta'_\infinity$, and so $\zeta_\infinity = \zeta'_\infinity
\in \Lambda$.

\subparagraph{Fact 4} $f_\infinity \from \bdy G \to \Lambda$ is a homeomorphism.

\medskip
Since both the domain and range are compact Hausdorff spaces it suffices to prove continuity in
one direction. Continuity of $f_\infinity^\inv$ follows by simply noting that for fixed $x \in \O$
and for a convergent sequence $\xi_n\to \xi$ in $\Lambda \subset \PMF$, the sequence of rays
$\ray{x}{\xi_n}$ converges in the compact open topology to the ray $\ray{x}{\xi}$.

\subparagraph{Fact 5} The map $\overline f = f \union f_\infinity \from G \union \bdy G \to
\overline\T = \T \union \PMF$ is continuous.

To be precise, this map is continuous using the Thurston compactification $\overline\T$ of $\T$.
We prove this by showing first that the map is continuous using the \Teichmuller\
compactification, and then we apply Masur's Two Boundaries Theorem \cite{Masur:TwoBoundaries}
which says that the map from the \Teichmuller\ compactification to the Thurston compactification
is continuous at uniquely ergodic points of $\PMF$.

First we recall the \Teichmuller\ compactification in a form convenient for our current purposes.
There are actually many different \Teichmuller\ compactifications, one for each choice of a base
point in $\T$; we shall fix a base point $z=f(x) \in \O$ for some $x \in G$. As we have seen,
there is a unique geodesic segment $[z,z']$ for each $z'\in\T$, and a unique geodesic ray
$\ray{z}{\zeta}$ for each $\zeta\in\PMF$. The \Teichmuller\ topology on $\overline\T = \T \union
\PMF$ restricts to the standard topologies on $\T$ and on
$\PMF$, it has $\T$ as a dense open subset, and a sequence $z_i\in \T$ converges to $\zeta \in
\PMF$ if and only if the sequence of segments $[z,z_i]$ converges to the ray $\ray{z}{\zeta}$ in
the compact open topology; equivalently, letting $B$ denote the unit ball in $\T$ centered on
$z$, the distance
$d(z,z_i)$ goes to infinity and the set $[z,z_i] \intersect B$ converges to the
set $\ray{z}{\zeta} \intersect B$ in the Hausdorff topology.

We already proved in Fact 4 that $f_\infinity$ is continuous; for this we implicitly used the fact
that the Thurston topology on $\PMF$ is identical to the \Teichmuller\ topology,
defined by identifying $\PMF$ with the unit tangent bundle at $x$. We also observed earlier, after
the proof that $f_\infinity$ is well-defined, that if $x_i \in G$ converges to $\xi \in \bdy G$,
then $f(x_i)\in\T$ converges to $f_\infinity(\xi) \in \PMF$ in the \Teichmuller\ topology on
$\overline\T$. Putting these together it follows that $\overline f$ is continuous using the
\Teichmuller\ topology on~$\overline\T$. Since $\Lambda=f_\infinity(\bdy G)$ consists entirely of
uniquely ergodic points in $\PMF$, Masur's Two Boundaries Theorem \cite{Masur:TwoBoundaries}
implies that the identity map on $\overline\T$ is continuous from the \Teichmuller\ topology to
the Thurston topology at each point of $\Lambda$, and so $\overline f$ is continuous using the
Thurston topology on~$\overline\T$.

\medskip

We now put the pieces together to complete the proof of convex cocompactness. Let $f' \from G \to
\WHull_G$ be an arbitrary $G$--equivariant map, and define $f_\infinity' \from \bdy G \to \PMF$ to be
equal to $f_\infinity$. We must prove that $f'$ is a quasi-isometry and that the extension $\bar f' =
f'\union f_\infinity' \from G \union\bdy G \to \WHull_G \union \Lambda_G$ is continuous.
From Facts 1--5 above, it follows that the quasi-isometry $f \from G
\to \union\G$ has continuous extension $\bar f \from G \union \bdy G \to \union\G \union
\Lambda$, and so $\union\G$ is a Gromov hyperbolic metric space with Gromov compactification
$\union\G \union \Lambda$. Since $\WHull_G \subset \union\G$ is a $G$--invariant subset, it
follows that $\WHull_G$ is Gromov hyperbolic with Gromov compactification $\WHull_G \union
\Lambda$. The map $f'$ is a $G$--equivariant map between quasigeodesic metric spaces on which $G$
acts properly and coboundedly by isometries, and hence $f'$ is a quasi-isometry. Since
$d(f'(x),f(x))$ is uniformly bounded for $x \in G$, then from the fact that $f_\infinity'=f_\infinity$ it
follows that $\bar f'$ is continuous. 

This completes the proof that weak orbit quasiconvexity implies convex cocompactness.

\paragraph{Convex cocompact implies weak orbit quasiconvexity} Assuming $G$ is convex
cocompact, pick a finite generating set for $G$ with Cayley graph $\Gamma$ and
$G$--equivariant, coarsely inverse quasi-isometries $f \from \Gamma \to \WHull_G$,
$\overline f \from \WHull_G \to \Gamma$.

Let $\O$ be an orbit of $G$ in $\T$. Since $G$ acts coboundedly on $\WHull_G$ it follows
that $\O$ has finite Hausdorff distance from $\WHull_G$ in $\T$. It suffices to show that for
any two points
$x,y \in \O$ there is a geodesic line whose infinite ends are in $\Lambda$ such that $x,y$
come within a uniformly finite distance of that line. 

Pick a $G$--equivariant map $g \from \Gamma \to \T$ taking the vertices of
$\Gamma$ bijectively to $\O$ and each edge of $\Gamma$ to a geodesic segment, so $f$ and $g$
differ by a bounded amount. Since $\Gamma$ is $\delta$--hyperbolic it follows that 
there is a constant $A$ such that any two vertices of $\Gamma$ lie within distance $A$ of some
bi-infinite geodesic. Pick $x,y \in \O$, and pick a bi-infinite geodesic $\gamma$ in $\Gamma$
such that $g^\inv(x), g^\inv(y)$ are within distance $A$ of $\gamma$. Let
$\xi,\eta \in \bdy G$ be the two ends of $\gamma$. By the statement of convex cocompactness,
there is a $K,C$ quasigeodesic line in $\Gamma$ of the form $\overline f \left(
\geodesic{f_\infinity\xi}{f_\infinity \eta} \right)$ whose two infinite ends are $\xi, \eta$,
where $K,C$ are independent of $\xi,\eta$. It follows that $\gamma$ and $\overline f \left(
\geodesic{f_\infinity\xi}{f_\infinity \eta} \right)$ are uniformly close, and so $f(\gamma)$
and $\geodesic{f_\infinity\xi}{f_\infinity \eta}$ are uniformly close, and so the points $x,
y$ are uniformly close to $\geodesic{f_\infinity\xi}{f_\infinity\eta}$.
\end{proof}

\section{Hyperbolic surface bundles over graphs}
\label{SectionBundlesOverGraphs}

In this section our goal is to give an explicit construction of model
geometries for surface group extensions, and to study regularity
properties of these geometries. Here is a brief outline; detailed
constructions follow.

Consider a finitely generated group $G$ and a homomorphism $f \from G
\to\Isom(\T) \approx \MCG$.  Let $X$ be a Cayley graph for $G$. Choose a
map $\Phi \from X \to \T$ which is equivariant with respect to the
homomorphism $f$, that is, $\Phi(g \cdot x) = f(g) \cdot \Phi(x)$, $x
\in X, g
\in G$, where we use the $\cdot$ notation to denote an action. By 
pulling back the canonical
marked hyperbolic surface bundle $\S\to\T$ via the map $\Phi$ we obtain
a marked hyperbolic surface bundle $\S_X \to X$. By pulling back the
canonical hyperbolic plane bundle $\H\to\T$ we obtain a hyperbolic plane
bundle $\H_X \to X$, and a covering map $\H_X\to \S_X$ with deck
transformation group $\pi_1(S)$. There is an action of the extension
group $\Gamma_G$ on $\H_X$ such that the covering map $\H_X \to \S_X$ is
equivariant with respect to the homomorphism $\Gamma_G \to G$.

By imposing a $G$--equivariant, proper, geodesic metric on $\S_X$ and lifting to $\H_X$, we can
then use $\H_X$ as a model geometry for the extension group~$\Gamma_G$. 

We may summarize all this in the following commutative diagrams:
$$\xymatrix{
          & \H \ar[r] \ar[dr] & \S \ar[d]  & & & \MCG(S,p) \ar[r] \ar[dr] & \MCG(S) \ar[d]\\
\H_X \ar[ur] \ar[r] \ar[dr] & \S_X \ar[ur] \ar[d] & \T  & & \Gamma_G
\ar[ur] \ar[r] \ar[dr] & G \ar[ur]_(.25)f \ar[d] & \MCG(S) \\
          & X \ar[ur]_\Phi & & & & G \ar[ur]_(.25)f
}$$
Each group in the right hand diagram acts on the corresponding space in the left hand diagram,
and each map in the left hand diagram is equivariant with respect to the corresponding group
homomorphism in the right hand diagram.

We will impose several $\Gamma_G$--equivariant structures on the space
$\H_X$, by finding appropriate $G$--equivariant structures on $\S_X$ and
lifting.
 
For example, we put an equivariant, proper, geodesic metric on $\H_X$ by
lifting an equivariant, proper, geodesic metric on $\S_X$. These metrics
will have the property that the topological fibrations $\S_X\to X$,
$\H_X \to X$ are also ``metric fibrations'' in the following sense. In a
metric space $Z$, given subsets $A,B \subset Z$, denote the \emph{min
distance} by $d_{\min}(A,B) =\inf\{d(a,b) \suchthat a \in A, b \in B\}$,
and the \emph{Hausdorff distance} by $d_\Haus(A,B) =\inf\{r\suchthat A
\subset N_r(B), B \subset N_r(A)\}$.

\begin{description}
\item[Metric fibration property] A map of metric spaces $f \from Z \to Y$ satisfies the
\emph{metric fibration property} if $Y$ is covered by neighborhoods $U$ such that if $y,z \in U$
then
$$d_{\min}(f^\inv(y),f^\inv(z)) = d_\Haus(f^\inv(y),f^\inv(z)) = d_Y(y,z)
$$
\end{description}

\subsection{Metrics and connections on surface bundles over paths}
\label{SectionBundlesOverArcs}

\paragraph{The marked hyperbolic surface bundle over a path in $\T$}

Consider first a smooth path $\alpha \from I \to \T$, defined on a
closed connected subset $I
\subset\reals$, that is, a closed interval, a closed ray, or the whole line. Pulling back the
canonical marked hyperbolic surface bundle $\S\to\T$ via the map $\alpha$ we obtain a marked
hyperbolic surface bundle $\S_\alpha\to I$. We impose a Riemannian metric on $\S_\alpha$ as
follows. 

Recall that we have chosen a connection on the bundle $\S\to\T$. By
pulling back the connection on the bundle $\S\to\T$ we obtain a
connection on the bundle $\S_\alpha\to I$, that is, a
\nb{1}dimensional sub-bundle of $T\S_\alpha$ which is complementary to the vertical sub-bundle
$T_v\S_\alpha$. There is a unique vector field $V$ on $S_\alpha$ parallel to the connection such
that the projection map $S_\alpha \to I$ takes each  vector of $V$ to a positive unit vector in
the tangent bundle of $I \subset \reals$. There is now a unique Riemannian metric on $\S$ whose
restriction to $T_v\S_\alpha$ is the given hyperbolic metric along leaves of $\S_\alpha$, and such
that $V$ is a unit vector field orthogonal to $T_v\S_\alpha$. Since $I$ is closed subset of
$\reals$, the path metric on $\S_\alpha$ induced from this Riemannian metric is proper, and so by
Fact~\ref{FactGeodesicMetricSpace} we may regard $\S_\alpha$ as a geodesic metric space.

Here is another description of the Riemannian metric on $\S_\alpha$. Integration of the
connection sub-bundle defines a \nb{1}dimensional foliation on $\S_\alpha$ transverse to the
surface fibration, whose leaves are called \emph{connection paths}. Choosing a base leaf of the
fibration $\S_\alpha \to I$, and identifying this base leaf with $S$, we may project along
connection paths to define a fibration $\S_\alpha \to S$. Combining this with the fibration
$\S_\alpha \to I$ we obtain a diffeomorphism $\S_\alpha \approx S \cross I$. Letting $g_t$ be the
given Riemannian metric of curvature $-1$ on the leaf $\S_t \approx S \cross t$, $t \in I$, we
obtain the Riemannian metric on $\S_\alpha$ via the formula
$$ds^2 = g_t^2 + dt^2.
$$

\subparagraph{Remark} The metric on $\S_\alpha$ depends on the choice of a connection on the
bundle $\S\to \T$. However, when $\alpha$ is cobounded, two different connections on $\S\to\T$
will induce metrics on $\S_\alpha$ which are bilipschitz equivalent, with bilipschitz constant
depending only on the pair of connections and on the coboundedness of $\alpha$, not on $\alpha$
itself.
\medskip

For each $s,t \in I$ we have a \emph{connection map} $h_{st} \from\S_s \to \S_t$, defined by
moving each point of $\S_s$ along a connection path until it hits $\S_t$. Clearly we have
$h_{st} \composed h_{rs} = h_{rt}$, ($r,s,t \in I$). Notice that the map $h_{st}$ takes each point
of $\S_s$ to the unique closest point on $\S_t$, and that point is at distance $\abs{s-t}$. In
fact, starting from an arbitrary point on $\S_s$, all paths to $\S_t$ have length $\ge
\abs{s-t}$, and the connection path is the unique one with length $=\abs{s-t}$. It follows that
the map $\S_\alpha \to I$ satisfies the metric fibration property.

Consider more generally a piecewise smooth path $\alpha \from I \to \T$. On each subinterval $I'
\subset I$ over which $\alpha$ is smooth, there is a Riemannian metric as constructed above. At a
point $t\in I$ where two such subintervals meet, the Riemannian metrics on the two sides agree
when restricted to $\S_t$. We therefore have a piecewise Riemannian metric on $\S_\alpha$,
inducing a proper geodesic metric. The connection paths which are defined over
each smooth subinterval $I' \subset I$ piece together to give connection paths on all of
$\S_\alpha$, and we obtain connection maps $h_{st} \from \S_s \to \S_t$ for all $s,t \in I$.

Note that since the connection on $\S\to\T$ is equivariant with respect to the action of $\MCG$,
the piecewise Riemannian metric on each $\S_\alpha$ is \emph{natural}, meaning that for
any $h\in \MCG$, the induced map $\S_\alpha \to \S_{h\composed\alpha}$ is an isometry. Similarly,
the connection paths and connection maps are also natural.

Each connection map $h_{st} \from \S_s \to \S_t$ is clearly a diffeomorphism, and since its domain
is compact it follows that $h_{st}$ is bilipschitz. The next proposition exhibits some
regularity, bounding the bilipschitz constant of $h_{st}$ by a function of $\abs{s-t}$ that
depends only on the coboundedness of the path $\alpha \from I \to \T$, and a lipschitz constant
for $\alpha$. For technical reasons we state the lemma only for paths $\alpha \from I \to \T$
which are \emph{piecewise affine}, meaning that $I$ is a concatenation of subintervals $I'$ such
that $\alpha \restrict I'$ is an \emph{affine} path, a constant speed
reparameterization of a \Teichmuller\ geodesic. Piecewise affine paths are sufficient for all of
what follows.

\begin{lemma}
\label{LemmaConnection}
For each bounded subset $\B \subset \Mod$ and each $\rho \ge 1$ there exists $K \ge 1$ such that
the following happens. If $\alpha \from I \to \T$ is a \nb{$\B$}cobounded, \nb{$\rho$}lipschitz,
piecewise affine path, then for each $s,t \in I$ the connection map $h_{st} \from \S_s \to
\S_t$ is $K^{\abs{s-t}}$--bilipschitz.
\end{lemma}

In what follows we shall describe the conclusion of this proposition by saying that $K$ is a
\emph{bilipschitz constant} for the connection maps on $\S_\alpha$. 

\begin{proof} A standard lemma found in most O.D.E.\ textbooks shows that if $\Phi$ is a
smooth flow on a compact manifold then there is a constant $K \ge 1$ such that $\norm{\Phi_t(v)}
\le K^{\abs{t}} \norm{v}$. We can plug into this argument as follows.

The conclusion of the lemma is local, and so it suffices to prove it under the assumption that
$I=[0,1]$ and that $\alpha$ is affine. There exists a compact subset $\A \subset \T$ such that
any $\B$--cobounded, $\rho$--lipschitz path $\alpha \from [0,1] \to \T$, can be translated by the
action of $\MCG$ to lie in the set $\A$. Let $C(\A,\rho)$ be the set of all
$\rho$--lipschitz affine paths $[0,1] \mapsto \A$, a compact space in the compact open topology.
By naturality of the metric on $\S_\alpha$, it suffices to prove the lemma for $\alpha\in
C(\A,\rho)$. For each $\alpha \in C(\A,\rho)$ and each vector $\vec w$ tangent to a fiber
$\S_s$, $s \in [0,1]$, define:
$$l(\vec w) = \lim_{t\to 0} \frac{1}{t} \log\left( \frac{\norm{D h_{s,s+t}(\vec w)}}{\norm{\vec
w}} 
\right)
   = \frac{d}{dt} \Biggm |_{t=0} \log\left( \frac{\norm{D h_{s,s+t}(\vec w)}}{\norm{\vec w}} 
\right)
$$
Since $l(c \vec w) = l(\vec w)$ for $c\ne 0$, we may regard $l(\vec w)$ as a function defined on
the projective tangent bundle of $S$ crossed with $I$, a compact space. As $\vec w$ varies,
and as $\alpha$ varies over the compact space $C(\A,\rho)$, the function $l(\vec w)$ varies
continuously, and so by compactness $l(\vec w)$ has a finite upper bound $l$. Setting $K=e^l$,
it now follows by standard methods that $\norm{h_{s,s+t}(\vec w)} \le K^{\abs{t}} \norm{\vec w}$
when $\vec w$ is tangent to $\S_s$, and so $h_{s,s+t}$ is $K^{\abs{t}}$ bilipschitz.
\end{proof}

\paragraph{The hyperbolic plane bundle over a path in $\T$}

Letting $\alpha\from I \to \T$ be a piecewise affine path as above, by
pulling back the canonical hyperbolic plane bundle $\H\to\T$ we obtain a
bundle $\H_\alpha \to I$. Note that there is a universal covering map
$\H_\alpha \to \S_\alpha$ with deck transformation group $\pi_1(S)$ such
that the composition $\H_\alpha \to \S_\alpha \to \S$ equals the
composition $\H_\alpha \to \H
\to \S$, and also the composition $\H_\alpha \to \S_\alpha \to I$ equals the fibration map
$\H_\alpha \to I$. By lifting the piecewise Riemannian metric from
$\S_\alpha$ we obtain a piecewise Riemannian metric on $\H_\alpha$,
inducing a proper, geodesic metric. The map $\H_\alpha
\to I$ satisfies the metric fibration property. The connection paths on $\S_\alpha$ lift to
connection paths on $\H_\alpha$, and we obtain connection maps $h_{st}
\from \H_s \to \H_t$. By applying Lemma~\ref{LemmaConnection} it follows
that if $\alpha$ is $\B$--cobounded and $\rho$--lipschitz then the same
constant $K=K(\B,\rho)$ is a bilipschitz constant for the connection
maps on $\H_\alpha$.

\subsection{Metrics and connections on surface bundles over graphs} 
\label{SectionModelSpace}

Let $f \from G \to \MCG$ be a homomorphism defined on a finitely
generated group~$G$. We have a canonical extension $1 \to \pi_1(S) \to
\Gamma_G \to G \to 1$.

Fix once and for all a Cayley graph $X$ for $G$, on which $G$ acts
cocompactly with quotient a rose. Fix a geodesic metric on $X$ with each
edge having length~1. Choose a $G$--equivariant map $\Phi \from X\to \T$
taking each edge of $X$ to an affine path in $\T$. Letting $\norm{\Phi}$
be the maximum speed of the map $\Phi$, ie,\ the maximal length of the
image of an edge of $X$ under $\Phi$, it follows that $\Phi$ is a
$\norm{\Phi}$--lipschitz map. Evidently the image of $\Phi$ is a
cobounded subset of $\T$, because the vertices of $X$ map to a single
orbit and each edge of $X$ maps to a geodesic of length
$\le\norm{\Phi}$. Choose a compact set $\B \subset \M$ so that
$\image(\Phi)$ is $\B$--cobounded.

Using the method of Section~\ref{SectionBundlesOverArcs}, for each edge
$e$ of $X$ we have a bundle $\S_e \to e$ equipped with a Riemannian
metric. Given any vertex $v$ of $X$, for any two edges $e,e'$ incident
to $v$ the Riemannian metrics on $\S_e$ and $\S_{e'}$ fit together
isometrically at $\S_v$. We may therefore paste together the Riemannian
metrics on $\S_e$ for all edges $e$ to obtain a marked hyperbolic
surface bundle $\S_X \to X$ equipped with a piecewise Riemannian
metric. The induced path metric on $\S_X$ is a proper, geodesic
metric. By naturality of the metrics on the bundles $\S_e$, the action
of $G$ on $X$ lifts to an isometric action on $\S_X$.

By lifting the metric from $\S_X$ to its universal cover $\H_X$ we
obtain a hyperbolic plane bundle $\H_X \to X$ on which the extension
group $\Gamma_G$ acts cocompactly, 
equipped with a $\Gamma_G$ equivariant, piecewise
Riemannian metric, inducing a proper, geodesic metric on $\H_X$.  Note 
in particular that $\Gamma_G$ is thus quasi-isometric to $\H_X$.

Note that this construction produces bundles $\S_X\to X$ and $\H_X \to
X$ isomorphic to the pullback bundles described at the beginning of
Section~\ref{SectionBundlesOverGraphs}. Since each map $\S_e \to e$,
$\H_e \to e$ satisfies the metric fibration property, it follows that
the maps $\S_X \to X$, $\H_X \to X$ also satisfy that property.

The connections on the spaces $\S_e$, for edges $e$ of $X$, piece together to define a
$G$--equivariant connection on $\S_X$. To make sense out of this, we consider only the connection
map defined for a piecewise path $\gamma \from [a,b] \to X$, as follows. The bundle $\S_X \to X$
pulls back to give a bundle $\S_\gamma \to [a,b]$, and the connection paths over each edge of $X$
piece together to give connection paths on $\S_\gamma$, with an induced connection map $h_\gamma
\from \S_{\gamma(a)} \to \S_{\gamma(b)}$. It follows immediately from Lemma~\ref{LemmaConnection}
that $h_\gamma$ is $K^{\Length(\gamma)}$--bilipschitz, where $K=K(\B,\norm{\Phi})$.

By lifting to $\H_X$, for each piecewise geodesic path $\gamma \from [a,b] \to X$ we similarly
obtain a $K^{\Length(\gamma)}$ bilipschitz connection map $\wt h_{\gamma} \from \H_{\gamma(a)} \to
\H_{\gamma(b)}$.

\subsection{Large scale geometry of surface bundles over paths} 
\label{SectionBundlesOverLines}

Our goal now is to compare metrics on $\H_\gamma$ and $\H_\beta$ for paths $\gamma,
\beta$ in $\T$ which are closely related. 

Given a metric space $Z$, two paths $\gamma,\beta \from I \to Z$, and a constant $A \ge 0$, we
say that $\gamma,\beta$ are \emph{$A$--fellow travellers} if $d(\gamma(t),\beta(t)) \le A$ for
all $t \in I$. More generally, given paths $\gamma \from I \to Z$,
$\beta \from J \to Z$, a constant $A \ge 0$, and constants $\lambda \ge 1, \epsilon \ge 0$, we
say that $\gamma$, $\beta$ are \emph{asynchronous $A$--fellow travellers} with respect to a
$\lambda,\epsilon$ quasi-isometry $\phi \from I \to J$ if the paths $\gamma$ and
$\beta\composed\phi$ are $A$--fellow travellers. It is a well known and simple fact that given
a quasigeodesic $\gamma \from I \to Z$ and another path $\beta \from J \to Z$, the following
are equivalent: 
\begin{enumerate}
\item $\beta$ is a quasigeodesic and $\beta,\gamma$ have finite Hausdorff distance;
\item $\beta$ is an asynchronous fellow traveller of $\gamma$.
\end{enumerate}
Moreover, the constants are uniformly related: in $1 \implies 2$, there exist asynchronous
fellow traveller constants $A,\lambda,\epsilon$ depending only on the quasigeodesic constants for
$\beta$ and the Hausdorff distance of $\beta,\gamma$; in $2 \implies 1$, there exist quasigeodesic
constants for $\beta$ and a bound on the Hausdorff distance between $\beta$ and $\gamma$
depending only on the asynchronous fellow traveller constants. 

The following proposition says that if $\gamma \from I \to \T$, $\beta \from J \to \T$ are
asynchronous fellow travellers in $\T$, then there is a fiber preserving quasi-isometry
$\H_\gamma \to\H_\beta$. Moreover, if $\gamma$ is a geodesic, and if instead of $\H_\gamma$ we
use the singular \solv\ space $\H^\solv_\gamma$, then there is a fiber preserving quasi-isometry
$\H^\solv_\gamma \to \H_\beta$.

\begin{proposition} 
\label{PropFellowTravellers}
For each bounded subset $\B \subset \Mod$, and each $\rho \ge 1$, $\lambda \ge 1$, $\epsilon
\ge 0$, $A \ge 0$, $K\ge 1$, there exists $K' \ge 1$, $C' \ge 0$ such that the following hold.
Suppose that $\gamma \from I \to \T$, $\beta\from J \to \T$ are $\B$--cobounded
$\rho$--Lipschitz, piecewise affine paths in $\T$. Suppose also that $\gamma,\beta$ are
asynchronous
$A$--fellow travellers, with respect to a $\lambda,\epsilon$ quasi-isometry
$\phi\from I \to J$. Then:
\begin{enumerate}
\item There exists a commutative diagram
$$\xymatrix{
\S_\gamma \ar[r]^{\Phi} \ar[d]  & \S_\beta \ar[d] \\
I    \ar[r]^{\phi}  &  J
}$$
such that the top row preserves markings, and such that any lifted map $\wt\Phi \from
\H_\gamma\to\H_\beta$ is a $K',C'$ quasi-isometry. 
\item If $\gamma$ is a geodesic, then there exists a commutative diagram
$$\xymatrix{
\S^\solv_\gamma \ar[r]^{\Phi} \ar[d] & \S_\beta \ar[d] \\
I          \ar[r]^{\phi} & J
}$$
such that the top row preserves markings, and such that any lifted map $\wt\Phi
\from \H^\solv_\gamma
\to \H_\beta$ is a $K',C'$ quasi-isometry.
\end{enumerate}
\end{proposition}

One way to interpret item (1) of this proposition is that a cobounded, lipschitz path in
\Teichmuller\ space has a well-defined geometry associated to it: approximate the given path by a
piecewise affine path and take the associated hyperbolic plane bundle; the metric on that bundle
is well-defined up to quasi-isometry, independent of the approximation. A further
argument shows that the geometry is independent of the choice of an equivariant connection
on the bundle $\S\to\T$: any two equivariant connections are related in a uniformly bilipschitz
manner over any cobounded subset of $\T$.

\begin{proof} Both (1) and (2) are proved in the same manner using
Proposition~\ref{PropTeichQuasi}; we prove only (1). 

To smooth the notation in the proof we denote $t' = \phi(t)$, we let $\S_t$ denote the
fiber $\S_{\gamma(t)}$ of $\S_\gamma$, we let $\S'_{t'}$ denote the corresponding
fiber $\S_{\beta(\phi(t'))}$ of $\S_\beta$, etc.

To prove (1), by applying Proposition~\ref{PropTeichQuasi}(1) we choose for each $t \in
\reals$ a marked map $\Phi_t \from \S_{t} \to \S'_{t'}$ for which any lift $\wt\Phi_t \from
\H_{t} \to \H'_{t'}$ is a $K_1,C_1$ quasi-isometry, where the constants $K_1,C_1$ depend only
on $\B,A$. Since each $\Phi_t$ preserves markings we may choose the lifts $\wt\Phi_t$ so that
for any $s,t$ we have a commutative diagram of induced boundary maps:
$$\xymatrix{
\bdy \H_{s} \ar[r]^{\bdy\wt\Phi_s} \ar[d]_{\bdy \wt h_{st}}  
                    & \bdy \H'_{s'} \ar[d]^{\bdy \wt h'_{s't'}} \\
\bdy \H_{t} \ar[r]_{\bdy\wt\Phi_{t}}   
                    & \bdy \H'_{t'}
}$$
Applying Proposition~\ref{PropTeichQuasi}(2) it follows that if we strip off the $\bdy$
symbols from the above diagram, and if we choose $s,t$ so that $\abs{s-t} \le 1$, then we
obtain the following diagram, a coarsely commutative diagram in the sense that the two paths
around the diagram differ in the sup norm by a constant $C_2$ depending only on
$\B,\rho,\lambda,\epsilon,A,K$:
$$\xymatrix{
\H_{s} \ar[r]^{\wt\Phi_s} \ar[d]_{\wt h_{s,t}}  
                    & \H'_{s'} \ar[d]^{\wt h_{s't'}}
\\
\H_{t} \ar[r]_{\wt\Phi_{t}}   
                    & \H'_{t'}
}$$
Define $\wt\Phi \from \H_\gamma \to \H_\beta$ so that $\wt\Phi \restrict \H_{s} =
\wt\Phi_s$. To prove that $\wt\Phi$ is a quasi-isometry we need only show that if $x,y \in
\H_\gamma$ satisfy $d(x,y) \le 1$ then $d(\wt\Phi(x),\wt\Phi(y))$ is bounded by a constant
depending only on $\B,\rho,\lambda,\epsilon,A,K$, and then carry out the similar argument
with inverses.

Given $x,y \in \H_\gamma$ with $d(x,y) \le 1$, choose $s,t$ so that $x \in \H_s$, $y \in
\H_t$. By the metric fibration property we have $\abs{s-t} \le 1$. Changing notation if
necessary we may assume that $s \le t$. Let $\alpha$ be the geodesic in $\H_\gamma$
connecting $x$ and $y$, and by the metric fibration property note that $\alpha \subset
\H_{[s-1,t+1]}$. Consider the map $p \from \H_{[s-1,t+1]} \to \H_t$ whose restriction to
$\H_r$ is the connection map $\wt h_{rt}$; it follows that $p$ is bilipschitz with constant
$K^{t-s+2} \le K^3$. The distance in $\H_t$ between the point $p(x) =
h_{st}(x)$ and the point $y$ is therefore at most $K^3$. Mapping over to $\H_\beta$ we have 
\begin{align*}
\begin{split}
d(\wt\Phi(x),\wt\Phi(y)) &\le d\bigl(\wt\Phi(x),h_{s't'}(\wt\Phi(x))\bigr) +
d\bigl(h_{s't'}(\wt\Phi(x)),\wt\Phi(h_{st}(x))\bigr)  \\
& \qquad\qquad \qquad \qquad \qquad  + d\bigl(\wt\Phi(h_{st}(x)),\wt\Phi(y)\bigr)
\end{split} \\
   &\le \abs{s'-t'} + C_2 + (K_1 K^3 + C_1)
\end{align*}
and since $\abs{s'-t'} \le \lambda \abs{s-t} + \epsilon \le \lambda+\epsilon$, the proof is
done.
\end{proof}

\section{Hyperbolic extension implies convex cocompact quotient}
\label{SectionExtensionQuotient}

In this section we prove Theorem~\ref{TheoremHypExtQuotient}. 

Fix a homomorphism $f \from G \to \MCG$ defined on a finitely generated
group~$G$, and suppose that the extension group $\Gamma_G$ is word
hyperbolic. We must prove that $f$ has finite kernel and that $f(G)$ is
a convex cocompact subgroup of $\MCG$.

Fix a Cayley graph $X$ for $G$ and an $f$--equivariant map $\Phi\from X
\to G$ which is affine on edges of $X$. Choose a bounded subset $\B
\subset \M$ and a number $\rho\ge 1$ such that $\Phi$ is $\B$--cobounded
and $\rho$--lipschitz. We have a hyperbolic plane bundle $\H_X \to X$,
and an action of $\Gamma_G$ on $\H_X$, such that the fibration $\H_X \to
X$ is equivariant with respect to the homomorphism $\Gamma_G \to G$. We
also have a piecewise Riemannian metric for which $\H_X \to X$ satisfies
the metric fibration property. We also have a connection on $\H_X$, in
the form of a connection map $h_\gamma \from \H_{\gamma(a)} \to
\H_{\gamma(b)}$ for any geodesic path $\gamma \from [a,b] \to X$. The
connection and metric are each equivariant with respect to
$\Gamma_G$. Since $\H_X$ is a proper geodesic metric space, it follows
that $\H_X$ is a model geometry for $\Gamma_G$. Since $\Gamma_G$ is word
hyperbolic, it follows that $\H_X$ is $\delta$--hyperbolic for some
$\delta \ge 0$.

\begin{fact}
\label{FactUnifProp}
For each point $x \in X$, the inclusion map $\H_x \inject \H_X$ is uniformly proper, with uniform
properness data independent of $x$.
\end{fact}

\begin{proof}
This follows because the subgroup of $\Gamma_G$ stabilizing $\H_x$ is the normal subgroup
$\pi_1(S)$, and the inclusion map $\pi_1(S) \inject \Gamma_G$ is uniformly proper with respect to
word metrics, a fact that holds for any finitely generated subgroup of a finitely generated group.
\end{proof}

For each geodesic path $\gamma \from I \to X$, $I$ a closed, connected subset of $\reals$, we
obtain a piecewise affine path $\Phi\composed\gamma \from I \to\T$ and a hyperbolic plane bundle
$\H_\gamma \to I$, which can be regarded either as the pullback of the bundle $\H\to\T$ via
$\Phi\composed \gamma$, or as the restriction of the bundle $\H_X \to X$ to $\gamma$. In either
case, we obtain a piecewise Riemannian metric and connection on $\H_\gamma$, natural with
respect to the action of $\pi_1(S)$. The connection on $\H_\gamma$ has bilipschitz constant $K$
depending only on $\B$ and $\rho$, meaning that for any $s,t \in\reals$, the connection map
$h_{st}\from\H_s \to \H_t$ is $K^{\abs{s-t}}$--bilipschitz.

Here is an outline of the proof of Theorem~\ref{TheoremHypExtQuotient}.

Our main task will be to prove that for each geodesic path $\gamma \from I \to X$, the space
$\H_\gamma$ is a $\delta'$--hyperbolic metric space, for some constant $\delta'$ depending only on
$\B$, $\rho$, and $\delta$. Of course, when $I$ is a finite segment the space $\H_\gamma$ is
quasi-isometric to the hyperbolic plane and so $\H_\gamma$ is a hyperbolic metric space, but
uniformity of the hyperbolicity constant $\delta'$ is crucial. This is obtained using the
concept of \emph{flaring}, introduced by Bestvina and Feighn for their combination theorem
\cite{BestvinaFeighn:combination}, and further developed by Gersten in
\cite{Gersten:cohomological}. The combination theorem says, in an appropriate context, that
flaring implies hyperbolicity. Gersten's converse, proved in the same context, says that
hyperbolicity implies flaring. We shall give a new technique for proving the converse, which
applies in a much broader, ``higher-dimensional'' context, and using this technique we show that
since $\H_X$ is $\delta$--hyperbolic it follows that each $\H_\gamma$ satisfies flaring, with
uniformity of constants. Then we shall apply the Bestvina--Feighn combination theorem in its
original context to conclude that $\H_\gamma$ is $\delta'$--hyperbolic.

Next we will apply a result of Mosher \cite{Mosher:StableQuasigeodesics} which says that since
$\H_\gamma$ is hyperbolic, the path $\Phi\composed\gamma \from I \to \T$ is a quasigeodesic
which is Hausdorff close to a \Teichmuller\ geodesic, again with uniformity of constants. This
will quickly imply finiteness of the kernel of $f$. The collection of these \Teichmuller\
geodesics, one for each geodesic $\gamma$ in $X$, will be used to verify the orbit
quasiconvexity property for the group $f(G)$.

In what follows, a path $I\xrightarrow{\gamma} X$ will often be confused with the composed
path $I\xrightarrow{\gamma} X \xrightarrow{\Phi} \T$; the context should make the meaning clear.

\paragraph{Remark} The context of the Bestvina--Feighn combination theorem, and Gersten's
converse, is the following. Consider a finite graph of groups $\Gamma$,
with word hyperbolic vertex and edge groups, such that each
edge-to-vertex group injection is a quasi-isometric
embedding. Associated to this is the Bass--Serre tree $T$, and a graph of
spaces $X \to T$ on which $\pi_1\Gamma$ acts properly discontinuously
and cocompactly. For each path in the tree $T$, Bestvina--Feighn define a
flaring condition on the portion of $X$ lying over that path. The
combination theorem combined with Gersten's converse says that flaring
is satisfied uniformly over all paths in the Bass--Serre tree if and only
if $\pi_1\Gamma$ is word hyperbolic. When $G$ is a free group mapped to
$\MCG$ then the extension $1 \to \pi_1 S \to \Gamma_G \to G \to 1$ fits
into this context, because $\Gamma_G$ is the fundamental group of a
graph of groups with edge and vertex groups isomorphic to $\pi_1 S$, and
with isomorphic edge-to-vertex injections, where the underlying graph is
a rose with fundamental group $G$. This was the technique used in
\cite{Mosher:hypbyhyp} to construct examples where $\Gamma_G$ is word hyperbolic. When $G$ is not
free then this doesn't work, motivating our ``higher-dimensional'' version of Gersten's result.

\subsection{Flaring}
\label{SectionFlaring}

Motivated by the statement of the Bestvina--Feighn combination theorem, we make
the following definitions.

Consider a sequence of positive real numbers $(r_j)_{j \in J}$, indexed by a subinterval $J$
of $\Z$. 

The \emph{$L$--lipschitz condition} says that $r_i / r_j < L^{\abs{i-j}}$ for all $i,j$, or
equivalently $r_i / r_j < L$ whenever $\abs{i-j}=1$. 

Given $\kappa > 1$, an
integer $n \ge 1$, and $A \ge 0$, we say that $(r_j)$ satisfies the
\emph{$(\kappa,n,A)$--flaring property} if, whenever the three integers $j-n$,
$j$, $j+n$ are all in $J$, we have:
$$r_j > A \quad\text{$\implies$}\quad \Max\{r_{j-n},r_{j+n}\} \ge \kappa \cdot r_j
$$
The number $A$ is called the \emph{flaring threshold}. Having a positive flaring threshold $A$
allows the sequence to stay bounded by $A$ on arbitrarily long intervals. However, at any place
where the sequence has a value larger than $A$, exponential growth kicks in inexorably, in
either the positive or the negative direction. 

Consider a piecewise affine, cobounded, lipschitz path $\gamma \from I \to \T$ and the
corresponding hyperbolic plane bundle $\H_\gamma \to I$. A
\emph{$\lambda$--quasivertical path} in $\H_\gamma$ is a $\lambda$--lipschitz path $\alpha\from
I'\to\H_\gamma$, defined on a subinterval $I' \subset I$, which is a section of the projection map
$\H_\gamma\to I$. For example, a $\lambda$--quasivertical path is a connection path if and only if
it is $1$--quasivertical. Note that each $\lambda$--quasivertical path is a
$(\lambda,0)$--quasigeodesic. 

The \emph{vertical flaring} property for the fibration $\H_\gamma \to \gamma$ says that there
exists $\kappa > 1$, an integer $n \ge 1$, and a function $A(\lambda) \from [1,\infinity) \to
(0,\infinity)$, such that if $\alpha,\beta\from I\to \H_\gamma$ are two $\lambda$--quasivertical
paths with the same domain $I'$, then setting $J=I'\intersect\Z$ the sequence
$$d_j\bigl(\alpha(j),\beta(j)\bigr), \quad j \in J
$$
satisfies the $\kappa,n,A(\lambda)$ flaring property, where $d_j$ is the distance function on
$\H_j$, $j\in J$. One can check that if the vertical flaring property holds for some function
$A(\lambda)$ then it holds for a function which grows linearly.

\begin{lemma}[Hyperbolicity of $\H_X$ implies vertical flaring of $\H_\gamma$] 
\label{LemmaHypImpliesFlare}
\hfill\break
With notation as above, for every $\delta$ there exists $\kappa$, $n$, $A(\lambda)$ such that if
$\H_X$ is $\delta$--hyperbolic then for each bi-infinite geodesic $\gamma$ in $X$ the fibration
$\H_\gamma \to I$ satisfies $\kappa$, $n$, $A(\lambda)$ vertical flaring.
\end{lemma}

The intuition behind the proof is that the flaring property is exactly analogous to the geodesic
divergence property in hyperbolic groups, described by Cannon in \cite{Cannon:TheoryHyp}. The
geodesic divergence property says that in a $\delta$--hyperbolic metric space, if $p$ is a base
point and if $\alpha,\beta$ are a pair of geodesic rays based at $p$, and if $d_i$ is the shortest
length of a path between $\alpha(i)$ and $\beta(i)$ that stays outside of the ball of radius $i$
centered on $p$, then the sequence $d_i$ satisfies a flaring property with constants independent
of
$\alpha,\beta$. In our context, $\alpha$ and $\beta$ will no longer have one endpoint in
common. But the quasivertical property together with the metric fibration property give us just
what we need to adapt Cannon's proof of geodesic divergence given in \cite{Cannon:TheoryHyp},
substituting the geodesic triangles in Cannon's proof with geodesic rectangles.

\begin{proof} We use $d$ for the metric on $\H_X$.

First observe that any $\lambda$--quasivertical path $\alpha$ in $\H_\gamma$ is a
$(\lambda,0)$--quasigeodesic in $\H_X$, in fact
$$\abs{s-t} \le d(\alpha(s),\alpha(t)) \le \lambda \abs{s-t}
$$
The upper bound is just the fact that $\alpha$ is $\lambda$--lipschitz, and the lower bound follows
from the metric fibration property for $\H_X \to X$, together with the fact that $\gamma$ is a
geodesic in $X$.

Consider then a pair of $\lambda$--quasivertical paths $\alpha,\beta \from I' \to \H_\gamma$
defined on a subinterval $I' \subset I$, and let $J=I' \intersect Z = \{j_-,\ldots,j_+\}$. We
assume that $j_+-j_-$ is even and let $j_0 = \frac{j_+-j_-}{2} \in J$. For each $j \in J$ we have
a fiber $\H_j$ isometric to $\hyp^2$, with metric denoted $d_j$. We must prove that the sequence
$D_j = d_j(\alpha(j),\beta(j))$ satisfies $\kappa,n,A$ flaring, with $\kappa,n$ independent of
$\lambda$ and with $\kappa,n,A$ independent of 
$\alpha$, $\beta$, and $\gamma$.

For $j,k \in J$ let $h_{jk} \from \H_j \to \H_k$ be the connection map, a
$K^{\abs{j-k}}$ bilipschitz map. 

For each $j \in J$ we have an $\H_j$ geodesic $\rho_j \from [0,D_j] \to \H_j$ with endpoints
$\alpha(j)$, $\beta(j)$. 

\begin{claim} There is a family of quasivertical paths $v$ described as follows:
\begin{itemize}
\item For each $j \in J$ and each $t \in [0,D_j]$ the family contains
a unique quasivertical path $v_{jt}
\from [j_-,j_+] \to \H_\gamma$ that passes through the point $\rho_j(t)$. If we fix $j
\in J$, we thus obtain a parameterization of the family $v_{jt}$ by points $t \in [0,D_j]$.
\item The ordering of the family $v_{jt}$ induced by the order on $t \in [0,D_j]$ is
independent of $j$. The first path $v_{j0}$ in the family is identified with $\alpha$, and the
last path $v_{jD_j}$ is identified with $\beta$.
\item Each $v_{jt}$ is $\lambda'$--quasivertical, where $\lambda'$ depends only on $\lambda$ and
$K$.
\end{itemize}
\end{claim}

When $j$ is assumed fixed, we write $v_t$ for the path $v_{jt}$.

\begin{proof}[Proof of claim] Given $j-1,j \in J$, consider the following $(K,0)$--quasigeodesic in
$\H_j$:
$$\rho'_j = h_{j-1,j} \composed \rho_{j-1} \from [0,D_{j-1}] \to \H_j
$$ 
Since connection paths are geodesics, and since $\alpha,\beta$ are
$\lambda$--quasivertical, it follows that the endpoint $\rho'_j(0) = h_{j-1,j}(\alpha(j-1))$ and
the corresponding endpoint $\rho_j(0) = \alpha(j)$ have distance in $\H_X$ at most $\lambda+1$,
and similarly for the opposite endpoints $\rho'_j(D_{j-1}) = h_{j-1,j}(\beta(j-1))$ and
$\rho_j(D_j) = \beta(j)$. Each endpoint of $\rho'_j$ and the corresponding endpoint of $\rho_j$
therefore have distance in $\H_j$ bounded by a constant depending only on $\lambda$; this follows
from Fact~\ref{FactUnifProp}. Since the spaces $\H_j$ are all isometric to $\hyp^2$, it follows
that the Hausdorff distance between
$\rho_j$ and $\rho'_j$ in $\H_j$ is bounded by a constant depending only on $K$, $\lambda$, which
implies in turn that there is a quasi-isometric reparameterization $r_j \from [0,D_{j-1}] \to
[0,D_j]$ such that
$$d_j\left(\rho'_j(t),\rho_j(r_j(t)) \right) \le D
$$
where the constant $D$ and the quasi-isometry constants for $r_j$ depend only on $K$, $\lambda$.
By possibly increasing the quasi-isometry constants we may assume furthermore that $r_j$ is an
orientation preserving homeomorphism. It follows that we may connect the point $\rho_{j-1}(t)$ to
the point $\rho_j(r_j(t))$ by a $\lambda'$--quasivertical path defined over the interval $[j-1,j]
\subset \reals$, where $\lambda'$ depends only on $K$, $\lambda$; when $t=0$ we may choose the
path to be $\alpha \restrict [j-1,j]$, and when $t=D_{j-1}$ we may choose the path $\beta
\restrict [j-1,j]$. By piecing together these paths as $j$ varies over $J$, we obtain the
required family of paths $v$.
\end{proof}

We use $\delta$--hyperbolicity of $\H_X$ in the following manner. First, for any geodesic
rectangle $a*b*c*d$ in $\H_X$ it follows that any point on $a$ is within distance $2\delta$ of $b
\union c\union d$. Second, for any $(\lambda',0)$ quasigeodesic in $\H_X$, the Hausdorff distance
to any geodesic with the same endpoints is bounded by a constant $\delta_1$ depending only on
$\delta,\lambda'$. For any rectangle of the form $v*\sigma*w*\sigma'$ where $\sigma,\sigma'$ are
geodesics and $v,w$ are $(\lambda',0)$ quasigeodesics, it follows that any point on $v$ is within
distance $\delta_2 = 2\delta + 2\delta_1$ of $\sigma \union w \union \sigma'$. 

By Fact~\ref{FactUnifProp} there exists a constant $\delta_3$ such that:
$$\text{for all } j \in J, x,y \in \H_j, \text{ if } d(x,y) \le (1+\lambda')\delta_2 \text{ then }
d_j(x,y) \le \delta_3
$$

We are now ready to define the flaring parameters $\kappa,n,A$. Let
\begin{align*}
\kappa &= \frac{3}{2} \\
n      &= \left\lfloor  \delta_2 + 3\delta_3  \right\rfloor + 1  \\
A      &= \delta_3
\end{align*}
where $\lfloor x \rfloor$ is the greatest integer $\le x$. Assuming as we may that $j_{\pm} = j_0
\pm n$ (and so the Hausdorff distance between $\H_{j_0}$ and $\H_{j_\pm}$ in $\H_X$ equals $n$),
we must prove:
\begin{itemize}
\item if $D_{j_0} > A$ then $\max\{D_{j_-},D_{j_+}\} \ge \kappa D_{j_0}.$
\end{itemize}

\paragraph{Case 1\qua $\max\{D_{j_-},D_{j_+}\} \le 6\delta_3$} It
follows that there is a rectangle in $\H_X$ of the form $\alpha * \sigma_- * \beta * \sigma_+$
where $\sigma_\pm$ is a geodesic in $\H_X$ with the same endpoints as $\rho_{j_\pm}$, and where
$\sigma_\pm$ has length $\le 6\delta_3$. Consider now the point $\alpha(j_0)$, whose distance from
some point $z \in \sigma_- \union \beta \union \sigma_+$ is at most $\delta_2$. If $z \in
\sigma_-$ then it follows that 
$$d(\alpha(j_0),\H_{j_+}) \le \delta_2 + \frac{6\delta_3}{2} < n, $$ a
contradiction. We reach a similar contradiction if $z \in
\sigma_+$. Therefore $z \in \beta$. It follows that $z=\beta(s) \in
\H_s$ for some $s$ such that $\abs{s-j_0} \le \delta_2$, and so by
following along $\beta$ a length at most $\lambda'\delta_2$ we reach the
point $\beta(j_0)$.  This shows that $d(\alpha(j_0),\beta(j_0)) \le
(1+\lambda') \delta_2$, and so $D_{j_0} \le \delta_3$, that is, $D_{j_0}
\le A$.

\paragraph{Case 2\qua $\max\{D_{j_-},D_{j_+}\} \ge 3\delta_3$} In the family $v$,
we claim that there is a discrete subfamily $\alpha = v_{t_0}, v_{t_1},
\ldots, v_{t_K} = \beta$, with $t_0<t_1<\cdots<t_K$, such that the
following property is satisfied: for each $k=1,\ldots,K$, letting $$
\Delta_{k\pm} = d_{j_{\pm}} \left(v_{t_{k-1}}(j_\pm),v_{t_{k}}(j_\pm) \right) 
$$
then we have
$$
\max\{\Delta_{k-},\Delta_{k+}\} \in [3\delta_3,6\delta_3].
$$ By assumption of Case 2, the subfamily
$\{\alpha=v_{t_0},\beta=v_{t_1}\}$ has the property
$\max\{\Delta_{k-},\Delta_{k+}\} = \max\{D_{j_-},D_{j_+}\} \ge
3\delta_3$ (for $k=1$). Suppose by induction that we have a subfamily
$\alpha = v_{t_0}, v_{t_1}, \ldots, v_{t_K} = \beta$, with
$t_0<t_1<\cdots<t_K$, such that $\max\{\Delta_{k-},\Delta_{k+}\} \ge
3\delta_3$ for all $k$, but suppose that
$\max\{\Delta_{k-},\Delta_{k+}\} > 6\delta_3$ for some $k$. If, say,
$\Delta_{k_+} > 6\delta_3$, then we subdivide the geodesic segment
$\rho_{j_+}[v_{t_{k-1}}(j_+),v_{t_{k}}(j_+)]$ in half at a point $t\in
\rho_{j_+}$, yielding two subsegments of length $>3\delta_3$, and we add
the path $v_{j_+t}$ to our subfamily; similarly, if $\Delta_{k_-} >
6\delta_3$ then we subdivide the interval
$\rho_{j_-}[v_{t_{k-1}}(j_-),v_{t_{k}}(j_-)]$ in half. This process must
eventually stop, because $$K \le \frac{1}{3\delta_3} \left( D_{j_-} +
D_{j_+} \right) $$ thereby proving the claim.

From the exact same argument as in Case~1, using the fact that 
$$\max\{\Delta_{k-},\Delta_{k+}\} \le 6\delta_3,
$$  
it now follows that
$$\Delta_{k0} = d_{j_0} \left(  v_{t_{k-1}}(j_0),v_{t_{k}}(j_0)  \right) \le \delta_3
$$
for all $k=1,\ldots,K$. 

We therefore have:
\begin{align*}
D_{j_0} &= \sum_{k=1}^K \Delta_{k0} \le K \delta_3 \\
D_{j_-} + D_{j_+} &= \sum_{k=1}^K \Delta_{k-} + \Delta_{k+} \ge \sum_{k=1}^K \max\{ \Delta_{k-},
\Delta_{k+} \} \\
                  &\ge K \cdot 3\delta_3  \\
\max\{D_{j_-},D_{j_+}\} &\ge \frac{3}{2} K \delta_3 \\
                        &\ge \frac{3}{2} D_{j_0}
\end{align*}
This completes the proof of Lemma~\ref{LemmaHypImpliesFlare}.
\end{proof}

\paragraph{Remark} The argument given in Lemma~\ref{LemmaHypImpliesFlare}, while stated
explicitly only for groups of the form $\Gamma_G$, generalizes to a much broader context. Graphs
of groups, the context for the Bestvina--Feighn combination theorem
\cite{BestvinaFeighn:combination} and Gersten's converse \cite{Gersten:cohomological}, have been
generalized to triangles of groups by Gersten and Stallings \cite{Stallings:triangles}, and to
general complexes of groups by Haefliger \cite{Haefliger:orbihedra}. The arguments of
Lemma~\ref{LemmaHypImpliesFlare} will also apply to show that a developable complex of groups
with word hyperbolic fundamental group satisfies a flaring property over any geodesic in the
universal covering complex. A converse would also be
nice, giving a higher dimensional generalization of the Bestvina--Feighn combination theorem, but
we do not know how to prove such a converse, nor do we have any examples to which it might apply
(see Question~\ref{QuestionNonFreeConvexCocompact} in the introduction).

\medskip

Next we have:

\begin{lemma}[Flaring implies hyperbolic]
\label{LemmaFlaringHyperbolic}
For each bounded subset $\B \subset \Mod$, each $\rho \ge 1$, and each set of flaring data
$\kappa>1$, $n \ge 1$, $A(\lambda)$, there exists $\delta \ge 0$ such that the following holds. If
$\gamma\from I\to\T$ is a $\B$--cobounded, $\rho$--lipschitz, piecewise affine path defined on a
subinterval $I \subset \reals$, and if the metric fibration $\H_\gamma\to I$ satisfies $\kappa,
n, A(\lambda)$ vertical flaring, then $\H_\gamma$ is $\delta$--hyperbolic.
\end{lemma}

\begin{proof} This is basically an immediate application of the Bestvina--Feighn combination
theorem \cite{BestvinaFeighn:combination}. To be formally correct, some
remarks are needed to translate from our present geometric setting, of a
hyperbolic plane bundle $\H_\gamma \to I$, to the combinatorial setting of
\cite{BestvinaFeighn:combination}, and to justify that our vertical flaring property for
$\H_\gamma$ corresponds to the ``hallways flare condition'' of
\cite{BestvinaFeighn:combination}.

We may assume that the endpoints of the interval $I$, if any, are integers.

The first observation is that 
there is a $\pi_1(S)$--equivariant triangulation $\wt\tau$ of
$\H_\gamma$ with the following properties:
\begin{description}
\item[Graph of spaces] \quad
\begin{itemize}
\item For each $n \in J=I\intersect\Z$ 
there is a \nb{2}dimensional subcomplex $\wt\tau_n$ which is a
triangulation of the hyperbolic plane $\H_n$.
\item Each \nb{1}cell of $\wt\tau$ is either \emph{horizontal} (a \nb{1}cell of some
$\tau_n$), or \emph{vertical} (connecting a vertex of some $\tau_n$ to a vertex of some
$\tau_{n+1}$);
\item each \nb{2}cell 
of $\wt\tau$ is either horizontal (a \nb{2}cell of some $\wt\tau_n$), or
vertical (meaning that the boundary contains exactly two vertical
\nb{1}cells).
\end{itemize}
\item[Bounded combinatorics] There is an upper bound 
depending only on $\B$,
$\rho$ for the valence of each \nb{0}cell and the number of sides of
each \nb{2}cell.
\item[Quasi-isometry] The inclusion of the \nb{1}skeleton of $\wt\tau$ into $\H_\gamma$ is a
quasi-isometry with constants depending only on $\B$ and $\rho$.
\end{description}
To see why $\wt\tau$ exists as described, consider the marked hyperbolic surface bundle $\S_\gamma
\to I$. For each hyperbolic surface $\S_n$, $n\in J$, there is a geodesic triangulation
$\tau_n$ of $\S_n$ with one vertex, whose edges have length bounded only in terms of $\B$. It
follows that there are constants $K'$, $C'$ depending only on $\B$, such that if $\wt\tau_n$ is
the lifted triangulation in $\H_n$, then the inclusion of the \nb{1}skeleton of $\wt\tau_n$ into
$\H_n$ is a $(K',C')$ quasi-isometry. Then, regarding $\bigcup_{n\in J} \tau_n$ as a triangulation
of $\bigcup_{n\in J}\S_n$, we can extend to a cell-decomposition $\tau$ of $\S_\gamma$ which is a
graph of spaces of bounded combinatorics. The existence of $\tau$ uses the fact that each
connection map $h_{n,n+1} \from \S_n \to \S_{n+1}$ is $K$--bilipschitz, so by moving each vertex of
$\tau_n$ along a connection path into $\S_{n+1}$ and them moving a finite distance to a vertex of
$\tau_{n+1}$ we obtain a $(K'',C'')$--quasi-isometry $h'_{n,n+1} \from \wt\tau_n\to\wt\tau_{n+1}$,
with $(K'',C'')$ depending only on $K$, and from this we easily construct $\tau$ so that its lift
$\wt\tau$ has the desired properties.

The second observation is that vertical flaring in $\H_\gamma$ is equivalent to the ``hallway
flare condition'' of \cite{BestvinaFeighn:combination} for $\wt\tau$, and this equivalence is
uniform with respect to the parameters in each property. To see why, note that quasivertical paths
in $\H_\gamma$ correspond to \emph{thin paths} in $\wt\tau$ as defined implicitly in
\cite{BestvinaFeighn:combination} Section~2: an edge path $\alpha \from I'=[m,n] \to \wt\tau$
is $\rho$--thin if the restriction of $\alpha$ to each subinterval $[i,i+1]$ lies in
$\wt\tau_{[i,i+1]}$ and is a concatenation of at most $\rho$ edges. Under the quasi-isometry
$\wt\tau \to \H_\gamma$ and its coarse inverse $\H_\gamma \to \wt\tau$, $\lambda$--quasivertical
paths in $\H_\gamma$ correspond to $\rho$--thin paths with a uniform relation between $\lambda$ and
$\rho$. 

In order to complete the translation from the geometric setting to the
combinatorial setting, while the results of
\cite{BestvinaFeighn:combination} are stated only when $\wt\tau$ is the
universal cover of a finite graph of spaces, nevertheless, the proofs
hold as stated for any graph of spaces with uniformly bounded
combinatorics: \emph{all} the steps in the proof extend to such graphs
of spaces, regardless of the presence of a deck transformation group
with compact quotient. The conclusion of the combination theorem is the
$\delta'$--hyperbolicity of the
\nb{1}skeleton of $\wt\tau$, with $\delta'$ depending only on the flaring constants for
$\wt\tau$, which depend in turn only on $\B$, $\rho$, and the flaring
constants for $\H_\gamma$.  It follows that $\H_\gamma$ is $\delta$
hyperbolic with the correct dependency for the constant $\delta$.
\end{proof}

\subsection{Proof of Theorem~\ref{TheoremHypExtQuotient}}
\label{SectionExtQuotientProof}

We adopt the notation from the beginning of Section~\ref{SectionExtensionQuotient}: a
homomorphism $f \from G \to\MCG$ determining the group $\Gamma_G$, a Cayley graph $X$ for $G$,
and a piecewise affine $f$--equivariant map $\Phi \from X \to \T$ which is
$\B$--cobounded and $\rho$--lipschitz. We have already proved, in Section~\ref{SectionExtensions},
that word hyperbolicity of $\Gamma_G$ implies finiteness of the kernel of $f$.

Letting $X^0$ be the \nb{0}skeleton, on which $G$ acts transitiveily, it follows that $\Phi(X^0)$
is an orbit of $f(G)$ in $\T$. We prove that $f(G)$ is convex cocompact by proving that
$\Phi(X^0)$ satisfies orbit quasiconvexity.

Choose two points $x,y \in X^0$. Let $\gamma \from I \to X$ be a geodesic segment connecting $x$
to $y$. Consider the composed path $I \xrightarrow{\gamma} X \xrightarrow{\Phi} \T$, which
by abuse of notation we shall also denote $\gamma$. There is a corresponding hyperbolic
plane bundle $\H_\gamma \to I$. Recall that $\gamma$ is $\B$--cobounded and
$\rho$--lipschitz in $\T$, with $\B$, $\rho$ independent of $\gamma$. Now apply
Lemmas~\ref{LemmaHypImpliesFlare} and~\ref{LemmaFlaringHyperbolic}, to conclude that $\H_\gamma$
is $\delta$--hyperbolic, with $\delta$ independent of $\gamma$. 

Now we quote the following result to obtain a \Teichmuller\ geodesic:

\begin{theorem}{\rm\cite{Mosher:StableQuasigeodesics}}\qua  
\label{TheoremStableQuasigeodesics}
For every bounded set $\B \subset \M$, $\rho \ge 1$, and $\delta \ge 0$, there exists
$\lambda \ge 1$, $\epsilon > 0$, and $A$ such that the following hold. If $\gamma \from I
\to \T$ is $\B$--cobounded and $\rho$--lipschitz, and if $\H_\gamma$ is $\delta$--hyperbolic, then
$\gamma$ is a $(\lambda,\epsilon)$--quasigeodesic, and there exists a \Teichmuller\ geodesic $g$,
sharing any endpoints of $\gamma$, such that $\gamma$ and $g$ have Hausdorff distance at most $A$.
\qed\end{theorem}

Letting $g$ be the \Teichmuller\ geodesic connecting $x$ to $y$ provided by the theorem, it
follows that $g$ is contained in the $A+\rho$ neighborhood of $\Phi(X^0)$. Since $x,y \in
\Phi(X^0)$ are arbitrary, this proves orbit quasiconvexity, and so $f(G)$ is convex
cocompact.

\section{Schottky groups}
\label{SectionSchottky}

\begin{definition}
A \emph{Schottky subgroup} of $\MCG$ is a free, convex cocompact subgroup. 
\end{definition}
The limit set $\Lambda \subset \PMF$ of a Schottky subgroup is therefore a Cantor set, and
every nontrivial element is pseudo-Anosov.

In this section we prove Theorem~\ref{TheoremSchottkyHypExt}, that a surface-by-free group is
word hyperbolic if and only if the free group is Schottky. One direction is already proved by
Theorem~\ref{TheoremHypExtQuotient}, and so we need only prove that when $F \subset \MCG$ is a
Schottky subgroup then $\Gamma_F \approx \pi_1(S) \semidirect F$ is word hyperbolic. 

Continuing with earlier notation, let $\Lambda \subset \PMF$ be the limit set of $F$ with weak
hull $\WHull_\Lambda$. Let $\tree$ be a Cayley graph for the group $F$, a tree on which $F$
acts properly discontinuously with quotient a rose. Let $\Phi \from \tree \to \T$ be an
$F$--equivariant map, affine on each edge, and $\rho$--lipschitz for some $\rho \ge 1$. There is a
bounded subset $\B \subset \Mod$ so that both $\WHull_\Lambda$ and $\Phi(\tree)$ are
$\B$--cobounded. We have a hyperbolic plane bundle $\H_\tree \to \tree$, on which $\pi_1(S)
\semidirect F$ acts properly discontinuously and cocompactly, and we have a piecewise Riemannian
metric on $\H_\tree$ on which $\pi_1(S) \semidirect F$ acts by isometries. 

We must prove that $\H_\tree$ is $\delta$--hyperbolic. By the
Bestvina--Feighn combination theorem \cite{BestvinaFeighn:combination}, it
is enough to show that for each bi-infinite geodesic
$\gamma$ in $\tree$, the bundle $\H_\gamma \to \reals$ satisfies vertical flaring, with flaring
data $\kappa, n, A(\lambda)$ independent of the choice of $\gamma$ (see the proof of
Lemma~\ref{LemmaFlaringHyperbolic} for translating the combinatorial setting of
\cite{BestvinaFeighn:combination} to our present geometric setting). 

Since $F$ is convex cocompact, there is a geodesic line $g$ in $\WHull_\Lambda$ which has finite
Hausdorff distance from $\Phi(\gamma)$. Let $\H^\solv_g$ be the singular \solv--space thereby
obtained. By Proposition~\ref{PropFellowTravellers}, the closest point map $\gamma \to g$ lifts
to a quasi-isometry $\H_\gamma \to \H^\solv_g$, with quasi-isometry constants independent of
$\gamma$, depending only on $\B$ and $\rho$. It therefore suffices to check the flaring
condition in $\H^\solv_g$, with flaring data independent of anything.

Take any $\kappa$ with $1 < \kappa < \frac{e^2}{2\sqrt{2}}$, say
$\kappa=2.6$. Let $n=2$. We show that for any $\lambda$ there is
an $A$ such that any two $\lambda$ quasivertical lines in $\H^\solv_g$
satisfy the $(\kappa,2,A)$--flaring condition. For this argument we do
not need that $g$ is cobounded (although in that case $\H^\solv_g$ may
not have bounded geometry).

Let $\alpha,\alpha' \from [-2,2] \to\H^\solv_g$ be two $\lambda$ quasivertical lines, lying over
a length 4 subsegment $[r-2,r+2]$ of $g \homeo \reals$. Let $x_i,y_i$ be the points where
$\alpha,\alpha'$ respectively intersect $\H_{r+i}$. Let $\xi_0=x_0$ and let $\xi_i$ be obtained by
flowing $x_0$ vertically into $\H_{r+i}$; define $\eta_0=y_0$ and $\eta_i$ similarly. Note that
for $i \in [-2,2]$ the points $\xi_i$ and $x_i$ are connected in $\H^\solv_g$ by a path which goes
along $\alpha$ from $\xi_i$ to $\xi_0$ travelling a distance at most $2\lambda$, and then
vertically from $\xi_0=x_0$ to $x_i$; the vertical projection of this path into $\H_i$ has length
at most $2e^2\lambda$, and so $d_i(x_i,\xi_i) \le 2e^2\lambda$. Similarly, $d_i(y_i,\eta_i) \le
2e^2\lambda$. 
 
We turn for the moment to showing that the sequence
$$d_{r+i}(\xi_{i},\eta_{i}), \quad i=-2,-1,0,1,2
$$
satisfies the $(\frac{e^2}{2\sqrt{2}},2,0)$--flaring condition. In the singular Euclidean surface
$\H_{r+i}$, let $\ell_i$ be the geodesic from $\xi_i$ to $\eta_i$, so the above sequence
becomes:
$$\Length(\ell_i), \quad i=-2,-1,0,1,2
$$
The singular Euclidean geodesic $\ell_0$ is a concatenation of subsegments of constant slope, two
consecutive subsegments meeting at a singularity. If at least half of $\ell_0$ has slope of
absolute value $\ge 1$ then:
$$
\frac{1}{2} \Length(\ell_0) \cdot \frac{1}{\sqrt{2}} \cdot e^2  \le \Length(\ell_2)
$$
If at least half of $\ell_0$ has slope of absolute value $\le 1$, we get a similar inequality
but with $\Length(\ell_{-2})$ on the right hand side. We have therefore shown:
$$\max\{d_{r+2}(\xi_2,\eta_2),d_{r-2}(\xi_{-2},\eta_{-2})\} \ge \frac{e^2}{2\sqrt{2}}
d_0(\xi_0,\eta_0)
$$
It follows that
\begin{align*}
\max\{d_{r+2}(x_2,y_2),d_{r-2}(x_{-2},y_{-2})\} &\ge \frac{e^2}{2\sqrt{2}} d_0(x_0,y_0)
- 2e^2\lambda
\\ & \ge \kappa d_0(x_0,y_0)
\end{align*}
where the last inequality holds as long as:
$$d_0(x_0,y_0) \ge A = \frac{2e^2\lambda}{\frac{e^2}{2\sqrt{2}}-\kappa}
$$
This ends the proof that $\pi_1(S)\semidirect F$ is word
hyperbolic when $F$ is Schottky.

\section{Extending the theory to orbifolds}
\label{SectionOrbifolds}

In this section we sketch how the theory can be extended to \nb{2}dimensional orbifolds. We shall
consider only those compact orbifolds whose underlying \nb{2}manifold is closed, and whose
orbifold locus therefore consists only of cone points, what we shall call a \emph{cone orbifold}.
The reason for this restriction is that if the underlying \nb{2}manifold has nonempty boundary
then the orbifold does not support any pseudo-Anosov homeomorphisms, since the isotopy classes of
the boundary curves must be permuted.\footnote{While the monograph \cite{FLP} develops a kind of
pseudo-Anosov theory on a bounded surface, it is \emph{not} appropriate for our present
purposes.} 

As it turns out, the mapping class group and \Teichmuller\ space of a cone orbifold depend not on
the actual orders of the different cone points, but only on the partition of the set of cone
points into subsets of constant order. For example, a spherical orbifold with one $\Z/2$ cone
point and three $\Z/4$ cone points has the same mapping class group and \Teichmuller\ space as a
spherical orbifold with three $\Z/42$ cone points and one $\Z/1000$ cone point. The relevant
structures can therefore be described more directly and economically in the following manner.

Let $S$ be a closed surface, not necessarily orientable. Let $\P=\{P_i\}_{i \in I}$ be a
finite, pairwise disjoint collection of finite, nonempty subsets of $S$. Let $\Homeo(S,\P)$ be
the group of homeomorphisms of $S$ which leave invariant each of the sets $P_i$,
$i \in I$. Let $\Homeo_0(S,\P)$ be the component of the identity of $\Homeo(S,\P)$ with
respect to the compact open topology; equivalently, $\Homeo_0(S,\P)$ consists of all elements of
$\Homeo(S,\P)$ which are isotopic to the identity through elements of $\Homeo(S,\P)$. The mapping
class group is $\MCG(S,\P) = \Homeo(S,\P) / \Homeo_0(S,\P)$. 

To define the \Teichmuller\ space, first we must widen the concept of a conformal structure so
that it applies to non-orientable surfaces, and we do this by allowing overlap maps which are
anticonformal as well as conformal. The \Teichmuller\ space $\T(S,\P)$ is then defined to be the
set of conformal structures on $S$ modulo the action of $\Homeo_0(S,\P)$. Quadratic differentials
and measured foliations on $(S,\P)$ are defined using the usual local models at points of $S -
\union\P$, but at a point of $\P$ a quadratic differential can have the local model $z^{n-2}
dz^2$ for any $n \ge 1$; the horizontal measured foliation of $z^{n-2} dz^2$ is the local model
for an $n$--pronged singularity of a measured foliation. Thus, at a point of $\union\P$ a measured
foliation can have any number of prongs $\ge 1$, whereas a singularity in $S-\union\P$ must
have $\ge 3$ prongs as usual. With these definitions, \Teichmuller\ maps are defined as usual,
making $\T(S,\P)$ into a proper geodesic metric space on which $\MCG(S,\P)$ acts properly
discontinuously, but not cocompactly; also, pseudo-Anosov homeomorphisms of $(S,\P)$ are defined
as usual.

We shall assume that $(S,\P)$ actually supports a pseudo-Anosov homeomorphism which has an
$n$--pronged singularity with $n \ne 2$. This rules out a small number of special cases, as
follows. When $S$ is a sphere, $\union\P$ must have at least four points. When $S$ is a projective
plane, $\union\P$ must have at least two points. When $S$ is a torus or Klein bottle, $\union\P$
must have at least one point. When $S$ is the surface of Euler characteristic $-1$, namely the
connected sum of a torus and a projective plane, the curve along which the torus and the
projective plane are glued is actually a characteristic curve for $S$, meaning that it is
preserved up to isotopy by any mapping class; therefore, in order for $(S,\P)$ to support a
pseudo-Anosov homeomorphism, $\union\P$ must have at least one point.

Now we apply these concepts to \nb{2}dimensional cone orbifolds. Suppose $\O$ is a cone orbifold
with underlying surface $S$. Let $P_n$ be the set of $\Z/n$ cone points, and let
$\P=\{P_n\}_{n \ge 2}$. Then we may define the mapping class group $\MCG(\O)$ to be $\MCG(S,\P)$,
and the \Teichmuller\ space $\T(\O)$ to be $\T(S,\P)$. Note that with the restrictions above on
the type of $(S,\P)$, the orbifold $\O$ has negative Euler characteristic. It follows that if
$\wt\O \to \O$ is the orbifold universal covering map, then for any conformal structure on $\O$
the lifted conformal structure is isomorphic to the Riemann disc. It follows that any conformal
structure on $\O$ can be uniquely uniformized to produce a hyperbolic structure, with a cone
angle of $2\pi/n$ at each $\Z/n$ cone point. 

At this stage we must confront the fact that the universal extension for surface groups, as
formulated in Section~\ref{SectionExtensions}, must be reformulated before it can be applied
to orbifolds. The Dehn--Nielsen--Baer--Epstein theorem is still true, as long as one uses orbifold
fundamental groups: if $p$ is a generic point of the cone orbifold $\O$, and if
$\pi_1(\O,p)$ is the orbifold fundamental group, then we have $\MCG(\O) \approx
\Out(\pi_1(\O,p))$. However, the ``once-punctured'' mapping class group $\MCG(\O,p)$ is
\emph{not} isomorphic to $\Aut(\pi_1(\O,p))$. For example, take a based simple loop
$\ell$ which bounds a disc whose interior contains a single $\Z/n$ cone point. In the
group $\pi_1(\O,p)$, the loop
$\ell$ represents an element of order $n$, and under the usual injection $\pi_1(\O,p) \inject
\Aut(\pi_1(\O,p))$ we obtain an element of order $n$. However, the element of $\MCG(\O,p)$
obtained by pushing $p$ around $\ell$ has infinite order in $\MCG(\O,p)$.

To repair this we need another group to take over the role of $\MCG(\O,p)$. Let $\wt\Homeo(\O)$
denote the group of homeomorphisms of $\wt\O$ which are lifts of homeomorphisms of $\O$, that is,
a homeomorphism
$\tilde f \from \wt\O \to \wt\O$ is in the group $\wt\Homeo(\O)$ if and only if there exists a
homeomorphism $f \from \O \to \O$ such that the following diagram commutes:
$$\xymatrix{
\wt\O \ar[r]^{\tilde f} \ar[d] &  \wt\O \ar[d] \\
\O  \ar[r]^f               & \O
}$$
With respect to the compact open topology, $\wt\Homeo(\O)$ becomes a topological group.
Let $\wt\Homeo_0(\O)$ be the component of the identity
$\wt\Homeo(\O)$. Equivalently,
$\wt\Homeo_0(\O)$ is the subgroup of elements of $\wt\Homeo(\O)$ isotopic to the identity through
elements of $\wt\Homeo(\O)$; alternatively it is the subgroup of $\wt\Homeo(\O)$ acting trivially
on the circle at infinity of $\wt\O \approx \hyp^2$. Define 
$$\wt\MCG(\O) = \wt\Homeo(\O)/\wt\Homeo_0(\O).
$$
Note that universal covering map $\wt\O \to \O$ induces a surjective
homomorphism $\wt\MCG(\O) \to
\MCG(\O)$, and the kernel is the group of deck transformations, isomorphic to $\pi_1(\O)$. We now
have a natural isomorphism of short exact sequences
$$\xymatrix{
1 \ar[r] & \pi_1(\O) \ar[r] \ar@{=}[d] & \wt\MCG(\O) \ar[r] \ar@2{~>}[d] &
\MCG(\O)
\ar[r] \ar@2{~>}[d] & 1
\\ 1 \ar[r] & \pi_1(\O) \ar[r]        & \Aut(\pi_1(\O)) \ar[r]   & \Out(\pi_1(\O)) \ar[r] &
1 }
$$
where we have suppressed the generic base point needed to define $\pi_1(\O)$.

We are now in a position to state that our main results,
Theorem~\ref{TheoremQuasiconvex}, \ref{TheoremHypExtQuotient}, \ref{TheoremSchottkyHypExt}, and
\ref{TheoremSchottkyAbundance}, are true with the orbifold $\O$ in place of the surface $S$, and
the proofs are unchanged. Although the references that we quote are stated solely in terms of
surfaces, namely \cite{Minsky:quasiprojections} and \cite{MasurMinsky:complex1} for
Theorem~\ref{TheoremQuasiconvex}, \cite{Mosher:HypExt} for Theorem~\ref{TheoremHypExtQuotient},
and \cite{Mosher:hypbyhyp} for Theorem~\ref{TheoremSchottkyAbundance}, nevertheless all the
proofs in those references work just as well for orbifolds instead of surfaces.

\end{document}

%% file: gtspec.tex
\def\ifplaintex{\expandafter\ifx\csname documentclass\endcsname\relax}

\expandafter\ifx\csname epsfbox\endcsname\relax\input epsf\fi

\ifplaintex 
\else
\fi

\def\gt{{\mathsurround=0pt\it $\cal G\mskip-2mu$eometry \&\ 
$\cal T\!\!$opology}}        

\def\gtp{{\mathsurround=0pt\it $\cal G\mskip-2mu$eometry \&\ 
$\cal T\!\!$opology $\cal P\!$ublications}}  

\def\lognumber#1{\def\thelognumber{#1}}
\def\volumenumber#1{\def\thevolumenumber{#1}}
\def\papernumber#1{\def\thepapernumber{#1}}
\def\volumeyear#1{\def\thevolumeyear{#1}}

\def\pagenumbers#1#2{\def\startpage{#1}\def\finishpage{#2}}
\def\published#1{\def\publishdate{#1}}
\def\proposed#1{\def\theproposer{#1}}
\def\seconded#1{\def\theseconders{#1}}
\def\received#1{\def\receiveddate{#1}}

\def\accepted#1{\def\accepteddate{#1}}

\def\asciiaddress#1{\def\theasciiaddress{#1}}
\def\asciiemail#1{\def\theasciiemail{#1}}
\long\def\asciiabstract#1{\long\def\theasciiabstract{#1}}

\let\\\par\let\thevolumenumber\relax\let\thepapernumber\relax
\let\thevolumeyear\relax\let\thesamplenumber\relax\let\startpage\relax
\let\finishpage\relax\let\publishdate\relax\let\receiveddate\relax
\let\reviseddate\relax\let\accepteddate\relax\let\theasciititle\relax
\let\theasciiauthors\relax\let\theasciiaddress\relax
\let\theasciiabstract\relax
\let\theasciiemail\relax\let\theshortauthors\relax\let\theshorttitle\relax

\long\def\maketitlep{   

\count0=\startpage

\gt\hfill      
\hbox to 77pt{\vbox to 0pt{\vglue -15pt\epsfbox{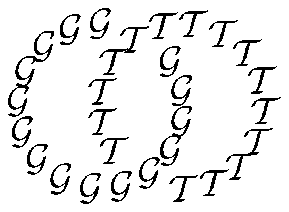}\vss}\hss}
\break
{\small\ifx\thesamplenumber\relax 
Volume \else Sample
\fi\thevolumenumber\ (\thevolumeyear)
\startpage--\finishpage\nl
Published: \publishdate}
\vglue 0.5truein plus 0.4fil minus 0.1truein

{\parskip=0pt\leftskip 0pt plus 1fil\def\\{\par\smallskip}{\ifplaintex\large
\else\Large\fi\bf\thetitle}\par\medskip}   

\vglue 0pt plus 0.1fil 

{\parskip=0pt\leftskip 0pt plus 1fil\def\\{\par}{\sc\theauthors}
\par\medskip}

\vglue 0pt plus 0.1fil 

{\small\parskip=0pt\let\newline\\
{\leftskip 0pt plus 1fil\def\\{\par}{\sl\theaddress}\par}
\expandafter\ifx\theemail\relax    
\relax\else\vglue 5pt plus 0.02fil minus 2pt\def\\{\stdspace{\rm 
and}\stdspace} 
\cl{Email:\stdspace\tt\theemail}\fi
\ifx\theurl\relax                  
\relax\else\vglue 5pt plus 0.02fil minus 2pt\def\\{\stdspace{\rm 
and}\stdspace}
\cl{URL:\stdspace\tt\theurl}\fi\par}

\vglue 7pt plus 0.3fil minus 3pt

{\bf Abstract}
\vglue 5pt plus 0.1fil minus 2pt

\theabstract

\vglue 7pt plus 0.3fil minus 3pt

{\bf AMS Classification numbers}\quad Primary:\quad \theprimaryclass

Secondary:\quad \thesecondaryclass

\vglue 5pt plus 0.3fil minus 2pt

{\bf Keywords:}\quad \thekeywords

\vglue 10pt plus 0.5fil minus 5pt

{\small  Proposed: \theproposer\hfill Received: \receiveddate\nl
Seconded: \theseconders\hfill 
\ifx\reviseddate\relax                         
Accepted: \accepteddate                        
\else
Revised: \reviseddate                          
\fi}
\eject
}       

\let\maketitlepage\maketitlep
\let\maketitle\maketitlepage

\font\phead=cmsl9 scaled 950
\font=cmsl9 scaled 1050
\font\pnum=cmbx10 scaled 913
\font\lnum=cmbx10 
\font\pfoot=cmsl9 scaled 950
\font=cmsl9 scaled 1050
\ifplaintex
\headline{\vbox to 0pt{\vskip -4.5mm\line{\small\phead\ifnum
\count0=\startpage ISSN 1364-0380 (on line)
1465-3060 (printed) \hfill {\pnum\folio}\else\ifodd\count0\def\\{ }%
\ifx\theshorttitle\relax\thetitle\else\theshorttitle\fi\hfill{\pnum\folio}
\else\def\\{ and }{\pnum\folio}\hfill\ifx\theshortauthors\relax\theauthors
\else\theshortauthors\fi\fi\fi}\vss}}
\footline{\vbox to 0pt{\vglue 0mm\line{\small\pfoot\ifnum\count0=\startpage
\copyright\ \gtp\hfill\else
\gt, Volume \thevolumenumber\ (\thevolumeyear)\hfill\fi}\vss
}}
\else
\makeatletter
\def\@oddhead{{\small\lhead\ifnum\count0=\startpage ISSN 1364-0380 (on line)
1465-3060 (printed) \hfill {\lnum\number\count0}\else\ifodd\count0
\def\\{ }\ifx\theshorttitle\relax \thetitle \else\theshorttitle\fi\hfill
{\lnum\number\count0}\else\def\\{ and }{\lnum\number\count0}
\hfill\ifx\theshortauthors\relax 
\theauthors\else\theshortauthors\fi\fi\fi}}\def\@evenhead{@oddhead}
\def\@oddfoot{\small\lfoot\ifnum\count0=\startpage\copyright\ \gtp\hfill\else
\gt, Volume \thevolumenumber\ (\thevolumeyear)\hfill\fi}
\def\@evenfoot{@oddfoot}
\makeatother
\fi

\newwrite\gtoutfile
\long\gdef\makeheadfile{  
{\def\\{, }\def\s{ }
\immediate\openout\gtoutfile head.xxx
\immediate\write\gtoutfile{To: math@arxiv.org}
\immediate\write\gtoutfile{Subject: put OR rep NNNNN:pppp}
\immediate\write\gtoutfile{--text follows this line--}
\immediate\write\gtoutfile{Proxy-for: \ifx\theasciiauthors\relax
\theauthors\else\theasciiauthors\fi\s<\ifx\theasciiemail\relax\theemail\else\theasciiemail\fi>}
\immediate\write\gtoutfile{\noexpand\\}
\immediate\write\gtoutfile{Authors: \ifx\theasciiauthors\relax
\theauthors\else\theasciiauthors\fi}
{\def\\{ }\immediate\write\gtoutfile{Title: \ifx\theasciititle\relax
\thetitle\else\theasciititle\fi}}
\immediate\write\gtoutfile{Subj-class: GT or GR or SG or ...}
\immediate\write\gtoutfile{MSC-class: \theprimaryclass\ifx\thesecondaryclass\relax\else, \thesecondaryclass\fi}
\immediate\write\gtoutfile{Journal-ref: Geom. Topol. \thevolumenumber\s
(\thevolumeyear) \startpage-\finishpage}
\immediate\write\gtoutfile{Comments: Published in Geometry and Topology at}
\immediate\write\gtoutfile{    http://www.maths.warwick.ac.uk/gt/GTVol\thevolumenumber/paper\thepapernumber.abs.html}
\immediate\write\gtoutfile{\noexpand\\}
\immediate\write\gtoutfile{}
\ifx\theasciiabstract\relax
\immediate\write\gtoutfile{\theabstract}\else
\immediate\write\gtoutfile{\theasciiabstract}\fi
\immediate\write\gtoutfile{}
\immediate\write\gtoutfile{\noexpand\\}
\immediate\write\gtoutfile{}
\immediate\closeout\gtoutfile}}  

\def\maketitlepage{\maketitlep\makeheadfile}
\let\maketitle\maketitlepage